\documentclass[12pt]{amsart}
\usepackage{amsmath,amsthm,amssymb}
\usepackage{mathrsfs}

\input epsf
\usepackage[all]{xy}
\newif\ifindex
\indextrue
\ifindex  
\else
   \renewcommand{\index}[1]{}
\fi
\ifindex\makeindex\fi

\newcommand{\Bc}[9]{\bibitem{#1} {#2}, \emph{#3}, in: \textbf{#4} (#5), #6 #7, #8--#9.}
\newcommand{\Emph}[1]{\emph{#1}\index{#1}}
\newcommand{\arx}[1]{\texttt{http://arxiv.org/abs/#1}}
\newcommand{\ed}{

\ifindex\printindex\fi

\end{document}}

      \newenvironment{changemargin}[2]{\begin{list}{}{
         \setlength{\topsep}{0pt}\setlength{\leftmargin}{0pt}
         \setlength{\rightmargin}{0pt}
         \setlength{\listparindent}{\parindent}
         \setlength{\itemindent}{\parindent}
         \setlength{\parsep}{0pt plus 1pt}
         \addtolength{\leftmargin}{#1}\addtolength{\rightmargin}{#2}
         }\item }{\end{list}}

\newcommand{\nop}{$\times$}
\newcommand{\yup}{\checkmark}
\newcommand{\mbq}{\mb{?}}
\newcommand{\arrays}{{{\{0,1\}}^{\N\x\N}}}
\newcommand{\CG}{C_\Gamma}
\newcommand{\CT}{C_\Tau}
\newcommand{\CO}{C_\Omega}
\newcommand{\psin}{pseudo-intersection}
\newcommand{\SPMBul}{\textbf{$\mathcal{SPM}$ Bulletin}}
\newcommand{\inv}{^{-1}}
\newcommand{\Cantor}{{\{0,1\}^\N}}
\newcommand{\SMZ}{\mathrm{SMZ}}
\newcommand{\CH}{the Continuum Hypothesis}
\newcommand{\diam}{\op{diam}}

\newcommand{\myeps}[1]{{#1}}
\newcommand{\CDR}{\mathsf{CDR}}

\newcommand{\compactN}{\mathbb{N}\cup\{\infty\}}
\newcommand{\NcompactN}{{(\compactN)^\N}}
\newcommand{\Increasing}{\mathcal{Z}}
\newcommand{\gone}{\mathsf{G}_1}
\newcommand{\gfin}{\mathsf{G}_{fin}}
\newcommand{\cl}[1]{\overline{#1}}
\newcommand{\Union}{\bigcup}

\newcommand{\fo}{\mathfrak{od}}
\newcommand{\h}{{\mathfrak h}}

\renewcommand{\u}{{\mathfrak u}}

\newcommand{\seq}[1]{\{#1\}_{n\in\N}}
\newcommand{\sr}[2]{{\txt{$#1$\\$#2$}}}
\newcommand{\scrA}{\mathscr{A}}
\newcommand{\scrB}{\mathscr{B}}
\newcommand{\scrC}{\mathscr{C}}
\newcommand{\scrD}{\mathscr{D}}
\newcommand{\cJ}{\mathcal{J}}
\newcommand{\cA}{\mathcal{A}}
\newcommand{\B}{\mathcal{B}}

\newcommand{\M}{\mathcal{M}}
\newcommand{\cU}{\mathcal{U}}
\newcommand{\cV}{\mathcal{V}}
\newcommand{\mb}[1]{{\mbox{\textbf{#1}}}}
\newcommand{\smb}[1]{{\!\!\mb{#1}\!\!}}

\newcommand{\fb}[1]{\framebox{$#1$}}

\renewcommand{\O}{\Cal O}

\newcommand{\BG}{\B_\Gamma}
\newcommand{\BT}{\B_\mathrm{T}}
\newcommand{\BO}{\B_\Omega}
\renewcommand{\split}{\mathsf{Split}}
\newcommand{\sone}{{\sf S}_1}
\newcommand{\sfin}{{\sf S}_{fin}}
\newcommand{\ufin}{{\sf U}_{fin}}
\long\def\forget#1\forgotten{}

\newcommand{\fin}{{[\N]^{<\aleph_0}}}
\newcommand{\roth}{{[\N]^{\aleph_0}}}
\newcommand{\Impl}{\Rightarrow}

\newcommand{\Cal}{\mathcal}

\newcommand{\cF}{\mathcal{F}}

\newcommand{\N}{{\mathbb N}}

\newcommand{\NN}{\N^\N}

\newcommand{\Q}{\mathbb Q}
\newcommand{\R}{\mathbb R}

\newcommand{\Tau}{\mathrm{T}}

\renewcommand{\b}{{\mathfrak b}}
\renewcommand{\c}{{\mathfrak c}}
\renewcommand{\d}{{\mathfrak d}}

\newcommand{\g}{\mathfrak{g}}

\newcommand{\oo}{\infty}
\newcommand{\p}{{\mathfrak p}}
\newcommand{\s}{\mathfrak{s}}

\newcommand{\w}{\omega}
\newcommand{\x}{\times}

\newcommand{\nin}{\not\in}

\newcommand{\sbst}{\subseteq}
\newcommand{\spst}{\supseteq}
\newcommand{\sm}{\setminus}
\newcommand{\as}{\subseteq^*}
\renewcommand{\|}{\upharpoonright}
\renewcommand{\(}{\left(}
\renewcommand{\)}{\right)}

\renewcommand{\>}{\rangle}

\newcommand{\op}{\operatorname}

\newcommand{\cov}{\op{\mathsf{cov}}}

\newcommand{\non}{\op{\mathsf{non}}}

\long\def\note#1\endnote%
{\ifNote{\par\medskip\tt Note: #1}\medskip\par
                         \else{}
                         \fi}

\newcommand{\impl}{\to}
\renewcommand{\t}{\mathfrak{t}}

\newtheorem{thm}{Theorem}[section]
\newtheorem{cor}[thm]{Corollary}

\newtheorem{lem}[thm]{Lemma}
\newtheorem{prop}[thm]{Proposition}
\newtheorem{prob}[thm]{Problem}

\theoremstyle{remark}

\newtheorem{rem}[thm]{Remark}

\theoremstyle{definition}
\newtheorem{defn}[thm]{Definition}

\newcommand{\be}{\begin{enumerate}}
\newcommand{\ee}{\end{enumerate}}
\newcommand{\bi}{\begin{itemize}}
\newcommand{\ei}{\end{itemize}}

\newcommand{\bpf}{\begin{proof}}
\newcommand{\epf}{\end{proof}}

\title{Some new directions in infinite-combinatorial topology}
\author{Boaz Tsaban}
\address{Department of Mathematics,
Weizmann Institute of Science,
Rehovot 76100,
Israel}
\email{boaz.tsaban@weizmann.ac.il}
\urladdr{http://www.cs.biu.ac.il/\~{}tsaban}

\begin{document}

\begin{abstract}
We give a light introduction to selection principles in topology,
a young subfield of infinite-combinatorial topology.
Emphasis is put on the modern approach to the problems it deals with.
Recent results are described, and open problems are stated.
Some results which do not appear elsewhere are also included, with proofs.
\end{abstract}

\maketitle

\tableofcontents

\setcounter{section}{-1}
\section{Introduction}
The modern era of what we call \Emph{infinite combinatorial topology},
or \Emph{selection principles in mathematics} began with Scheepers'
paper \cite{coc1} and the subsequent work \cite{coc2}.
In these works, a unified framework was given that extends many
particular investigations carried on in the classical era.
The current paper aims to give the reader a taste of the field
by telling six stories, each shedding light on one specific
theme. The stories are short, but in a sense never ending, since
each of them poses several open problems, and more are expected to
arise when these are solved.

This is not intended to be a systematic exposition to the field,
not even when we limit our scope to selection principles in \emph{topology}.
For that see Scheepers' survey \cite{LecceSurvey} as well as Ko\v{c}inac's \cite{KocSurv}.
Rather, we describe themes and results with which we are familiar.
This implies the disadvantage that we are often quoting our own results,
which only form a tiny portion of the field.\footnote{To partially compensate for this,
the name of the present author is never explicitly mentioned in the
paper.}

Some open problems are presented here. For many more see \cite{futurespm}.
After reading this introduction,
the reader can proceed directly to some of the
works of the experts in the field (or in closely related fields), such as:%
\footnote{This \emph{very} incomplete list is ordered alphabetically.
We did not give references to works not explicitly mentioned in this paper,
but the reader can find some of these in the bibliographies of the given references,
most notably, in \cite{LecceSurvey, KocSurv}.
We should also comment that not all works of the authors are formulated using the systematic notation of
Scheepers which we use here, and probably some of the mentioned experts do not consider themselves as
working in the discussed field (but undoubtly, each of them made significant contributions to the field).}
Liljana Babinkostova, Taras Banakh, Tomek Bartoszy\'nski,
Lev Bukovsk\'y, Krzysztof Ciesielski,
David Fremlin, Fred Galvin,
Salvador Garc\'ia-Ferreira, Janos Gerlits, Cosimo Guido,
Istvan Juhasz, Ljubisa Ko\v{c}inac, Adam Krawczyk,
Henryk Michalewski, Arnold Miller, Zsigmond Nagy,
Andrzej Nowik, Janusz Pawlikowski, Roman Pol,
Ireneusz Rec\l{}aw, Miroslav Repick\'y, Masami Sakai,
Marion Scheepers, Lajos Soukup, Paul Szeptycki,
Stevo Todorcevic, Tomasz Weiss, Lyubomyr Zdomskyy,
and many others.

\subsubsection*{Notation}
In most cases, the notation and terminology
we use is Scheepers' modern one, and we do not pay special attention to
the historical predecessors of the notations we use.
The reader can use the index at the end of the paper to locate the definitions
he is missing.

By \Emph{set of reals} we usually mean a zero-dimensional, separable metrizable space,
though some of the results hold in more general situations.

\subsubsection*{Apology}
Some of the results might be miscredited or misquoted (or both).
Please let us know of any mistake you find and we will
correct it in the online version of this paper \cite{ictonline}.

\section{The Menger-Hurewicz conjectures}

\subsection{The Menger and Hurewicz properties}
In 1924, Menger \cite{MENGER} introduced the following
basis covering property for a metric space $X$:
\begin{quote}
For each basis $\B$ of $X$, there exists a sequence $\seq{B_n}$ in
$\B$ such that $\lim_{n\to\infty}\diam(B_n) = 0$ and $X=\Union_nB_n$.
\end{quote}
It is an amusing exercise to show that every compact, and even
$\sigma$-compact\index{$\sigma$-compact} space has this property.
Menger conjectured that this property characterizes $\sigma$-compactness.
In 1925, Hurewicz \cite{HURE25} introduced two properties of the following
prototype.
For collections $\scrA,\scrB$ of covers of a space $X$, define
\bi
\item[$\ufin(\scrA,\scrB)$\index{$\ufin(\scrA,\scrB)$}:]
For each sequence $\seq{\cU_n}$ of members of $\scrA$
\emph{which do not contain a finite subcover},
there exist finite (possibly empty) subsets $\cF_n\sbst\cU_n$, $n\in\N$,
such that $\seq{\cup\cF_n}\in\scrB$.
\ei

\begin{figure}[!htp]
\begin{center}
\myeps{\epsfysize=6 truecm {\epsfbox{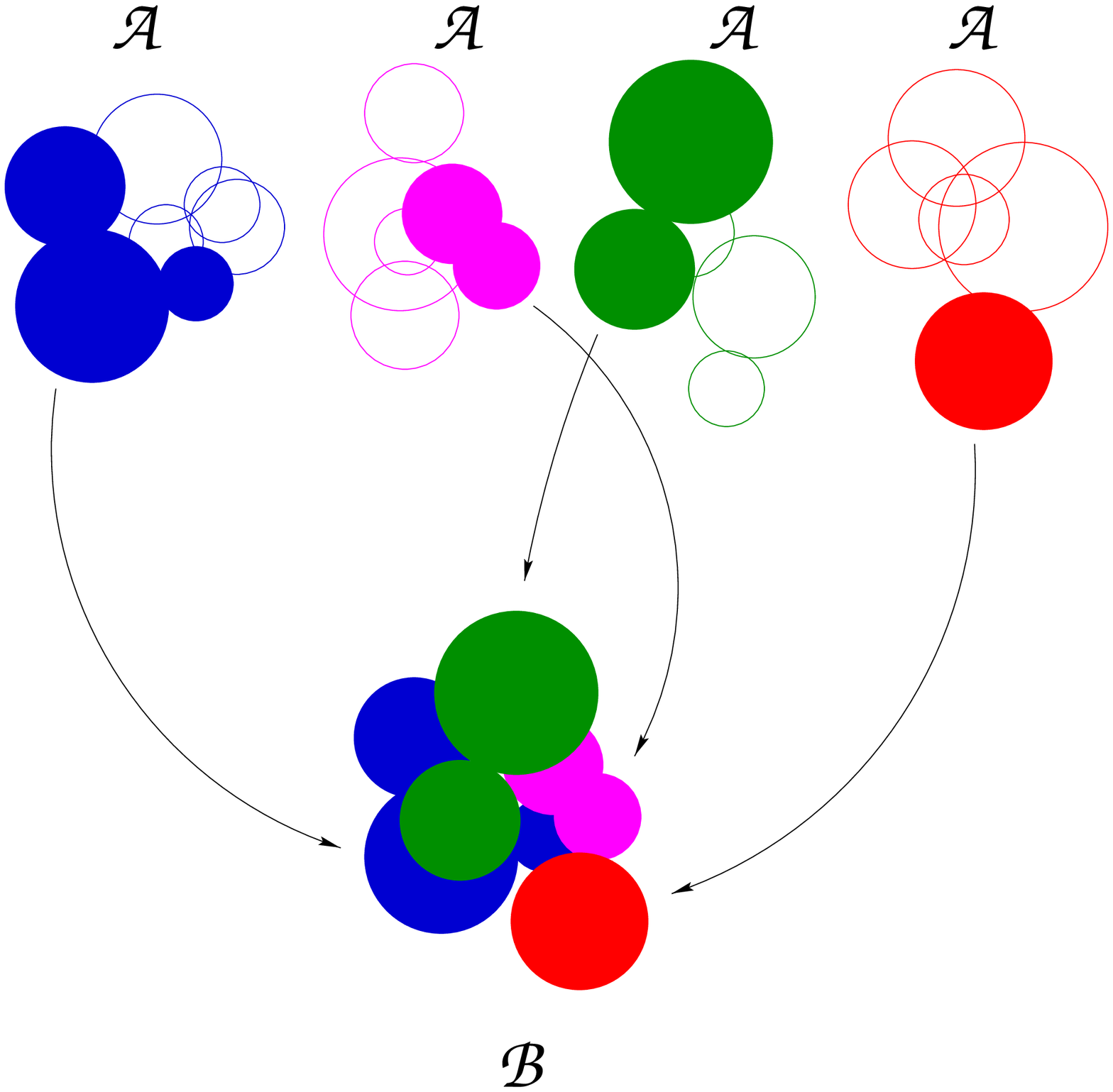}}}
\caption{$\ufin(\scrA,\scrB)$}
\end{center}
\end{figure}
Hurewicz proved that for $\O$\index{$\O$} the collection of \emph{all} open covers
of $X$, $\ufin(\O,\O)$ is equivalent to Menger's basis property.
Hurewicz did not settle Menger's conjecture, but he suggested a more modest one:
Call an open cover $\cU$ of $X$ a \Emph{$\gamma$-cover}
if $\cU$ is infinite, and each $x\in X$ belongs to all but finitely many
members of $\cU$.
Let $\Gamma$ denote the collection of all open $\gamma$-covers of $X$.
Clearly,
$$\index{$\sigma$-compact}\sigma\mbox{-compact}\Impl\ufin(\O,\Gamma)\Impl\ufin(\O,\O),$$
and Hurewicz conjectured that for metric spaces,
$\ufin(\O,\Gamma)$ (now known as the \Emph{Hurewicz property})
characterizes $\sigma$-compactness.

\subsection{Consistent counter examples}
It did not take long to find out that these conjectures are false
assuming \CH{}: Observe that every uncountable $F_\sigma$ set of reals contains an
uncountable perfect set, which in turn contains an uncountable set which is both
meager (i.e., of Baire first category) and null (i.e., of Lebesgue measure zero).
A set of reals $L$ is a \Emph{Luzin set} if it is uncountable, but for each meager set
$M$, $L\cap M$ is countable.
It was proved by Mahlo (1913) and Luzin (1914) that \CH{} implies the existence of a
Luzin set.
Sierpi\'nski pointed out that each
Luzin set has Menger's property $\ufin(\O,\O)$\label{LisMen}
(hint: If we cover a countable dense subset of the Luzin set by open sets, then the uncovered
part is meager and therefore countable),
and is therefore a counter example to Menger's conjecture.\footnote{%
This observation was ``added in proof'' just after footnote 1 on page 196
of Hurewicz' 1927 paper \cite{HURE27}.
}

Similarly, a set of reals $S$ is a \Emph{Sierpi\'nski  set} if it is uncountable,
but for each null set $N$, $S\cap N$ is countable.
Sierpi\'nski showed that \CH{} implies the existence of such sets,
and it can be shown that Sierpi\'nski sets have the Hurewicz property,
and is therefore a counter-examples to Hurewicz' (and therefore Menger's) conjecture.

Why must a Sierpi\'nski set satisfy Hurewicz' property $\ufin(\O,\Gamma)$?
A classical proof can be carried out using Egoroff's Theorem, but let us see
how a modern, combinatorial proof goes \cite{CBC}.
The \Emph{Baire space} $\NN$\index{$\NN$} (a Tychonoff power of the discrete space $\N$)
carries an interesting combinatorial structure:
For $f,g\in\NN$, write
$$f\le^* g\quad \mbox{if}\quad f(n)\le g(n)\mbox{ for all but finitely many }n.\index{$\le^*$}$$
$B\sbst\NN$ is \Emph{bounded} if
there exists $g\in\NN$ such that $f\le^* g$ for all $f\in B$.
$D\sbst\NN$ is \Emph{dominating} if
for each $g\in\NN$ there exists $f\in D$ such that $g\le^* f$.
It is easy to see that a countable union of bounded sets in $\NN$ is bounded,
and that compact (and therefore $\sigma$-compact\index{$\sigma$-compact}) subsets of $\NN$ are bounded.
The following theorem is essentially due to Hurewicz, who proved a variant of it in \cite{HURE27}.
In the form below, the theorem was stated and proved in Rec\l{}aw \cite{RECLAW}
in the zero-dimensional case, and then extended
by Zdomskyy \cite{ZdImages} to arbitrary subsets of $\R$.\footnote{%
In the current form, Theorem \ref{hure} does not hold for subsets of $\R^2$ \cite{ZdImages}.
}

\begin{thm}[Hurewicz]\label{hure}
For a set of reals $X$:
\be
\item $X$ satisfies $\ufin(\O,\Gamma)$ if,
and only if, all continuous images of $X$ in $\NN$ are bounded.
\item $X$ satisfies $\ufin(\O,\O)$ if,
and only if, all continuous images of $X$ in $\NN$ are not dominating.
\ee
\end{thm}

Assume that $S\sbst[0,1]$ is a Sierpi\'nski set
and $\Psi:S\to\NN$ is continuous. Then $\Psi$ can be extended to all of $[0,1]$
as a Borel function.
By a theorem of Luzin, there exists for each $n$ a closed subset $C_n$ of $[0,1]$
such that $\mu(C_n)\ge 1 - 1/n$,
and such that $\Psi$ is continuous on $C_n$. Since $C_n$ is compact,
$\Phi[C_n]$ is bounded in $\NN$.
The set $N=[0,1]\setminus \Union_n C_n$ is null, and so its intersection
with $S$ is countable. Consequently, $\Psi[S]$ is contained in a union of countably many bounded sets
in $\NN$, and is therefore bounded.

\subsection{Counter examples in ZFC}\label{ZFCexamples}
But are the conjectures \Emph{consistent}? It turns out that the answer is negative.
Surprisingly, it was only recently that this question was clarified, again using a combinatorial
approach.
Let $\b$\index{$\b$} denote the minimal size of an unbounded subset of $\NN$, and
$\d$\index{$\d$} denote the minimal size of a dominating subset of $\NN$.
The \Emph{critical cardinality} of a (nontrivial) collection $\cJ$ of sets of reals is
$$\non(\cJ)=\min\{|X| : X\sbst\R,\ X\nin\cJ\}.$$\index{$\non(\cJ)$}
By Hurewicz' Theorem \ref{hure},
$\non(\ufin(\O,\Gamma))=\b$, and $\non(\ufin(\O,\O))\allowbreak=\d$.

In 1988, Fremlin and Miller \cite{FM} used their celebrated dichotomic
argument\index{Fremlin-Miller dichotomy} to refute Menger's conjecture (in ZFC):
By the last observation,
if $\aleph_1<\d$ then any set of reals of size $\aleph_1$ will do;
and if $\aleph_1=\d$, then one can use a sophisticated combinatorial construction.
Chaber and Pol \cite{ChaPol}, exploiting a celebrated topological technique due to
Michael, extended the dichotomic argument to show that there always exists
a counter example of size $\b$ to Menger's conjecture.

Hurewicz' conjecture was refuted in 1996
\cite{coc2}, this time using the dichotomy $\aleph_1<\b$ or $\aleph_1=\b$ and
even more sophisticated combinatorial arguments in the second case.
This was improved by Scheepers \cite{wqn}, who showed (again on a dichotomic basis)
that there always exists a counter example of size $\t$\index{$\t$} ($\t$, to be defined in Section
\ref{MTP}, is an uncountable cardinal which is consistently greater than $\aleph_1$).

A simple construction was very recently found
\cite{BaShCon2000, ideals}
to refute both conjectures, and not on a dichotomic basis: There exists
a non $\sigma$-compact set of reals $H$ of size $\b$ which has
the Hurewicz property. The construction does not use any special hypothesis:
Let $\compactN$\index{$\compactN$} be the one point compactification of $\N$,
and $\Increasing$\index{$\Increasing$} be the nondecreasing functions $f\in\NcompactN$.
For a finite nondecreasing sequence $s$ of natural numbers, let
$q_{s}$ be the element of $\Increasing$ extending $s$ and being equal
to $\infty$ on all new $n$'s.
Then the collection $Q$ of all these elements $q_s$ is dense in $\Increasing$.
Define a
\emph{$\b$-scale}\index{$\b$-scale} to be an unbounded set $\{f_\alpha : \alpha<\b\}\sbst\NN$
of increasing functions, such that $f_\alpha\le^* f_\beta$ whenever $\alpha<\beta$.
It is an easy exercise to construct a $\b$-scale in ZFC.

\begin{thm}[Bartoszy\'nski, et.\ al.\ \cite{ideals}]\label{H}
Let $H$\index{Hurewicz set $H$} be a union of a $\b$-scale and $Q$.
Then all finite powers of $H$ satisfy $\ufin(\O,\Gamma)$
(but $H$ is not $\sigma$-compact).
\end{thm}

Consequently, there exists a counter example $G_H$\index{Hurewicz group $G_H$}
to the Hurewicz conjecture, such that $|G_H|=\b$ and $G_H$ is a
\emph{subgroup} of $\R$ \cite{o-bdd}.

Similarly, it was shown in \cite{ideals} that there exists a counter
example of size $\d$ to the Menger conjecture. However it is open
whether the group theoretic version also holds.

\begin{prob}[\cite{o-bdd}]\label{mengp}
Does there exist (in ZFC) a subgroup $G_M$\index{Menger group $G_M$}
of $\R$ such that $|G_M|=\d$ and $G_M$ has Menger's property $\ufin(\O,\O)$?
\end{prob}

\section{The Borel Conjecture}

\subsection{Strong measure zero}
Recall that a set of reals $X$ is \Emph{null}
if for each positive $\epsilon$ there exists
a cover $\seq{I_n}$ of $X$ such that $\sum_n\diam(I_n)<\epsilon$.
In his 1919 paper \cite{Borel}, Borel introduced the following
stronger property: A set of reals $X$ is \Emph{strongly null}
(or: has \Emph{strong measure zero}) if, for each sequence
$\seq{\epsilon_n}$ of positive reals, there exists a cover
$\seq{I_n}$ of $X$ such that $\diam(I_n)<\epsilon_n$ for all $n$.
But Borel was unable to construct a nontrivial (that is, an uncountable)
example of a strongly null set. He therefore conjectured that
there exist no such examples.
Sierpi\'nski (1928) observed that every Luzin set is strongly null
(see the hint on page \pageref{LisMen} for the reason),
thus \CH{} implies that \Emph{Borel's Conjecture} is false.
Sierpi\'nski asked whether the property of being strongly null
is preserved under homeomorphic (or even continuous) images.
The answer, given by Rothberger (1941) in \cite{ROTH41},
is negative under \CH{}.

If we carefully check Rothberger's argument,
we can obtain a slightly stronger result without making the proof
more complicated, and with the benefit of understanding the underlying
combinatorics better.
Theorem \ref{roth} and Proposition \ref{weiss} below are probably folklore, but we
do not know of a satisfactory reference for them so we give complete proofs.
(The proof of Theorem \ref{roth} is a modification of the proof of \cite[Theorem 2.9]{SMZjubar}.)
Let $\SMZ$\index{$\SMZ$} denote the collection of strongly null sets of reals.
A subset $A$ of $\NN$ is \Emph{strongly unbounded}
if for each $f\in\NN$, $|\{g\in A : g\le^* f\}|<|A|$.
Observe that there exist strongly unbounded sets of sizes $\b$ and $\d$,
thus any of the hypotheses $\non(\SMZ)=\b$ or $\non(\SMZ)=\d$
(in particular, \CH{}) implies the assumption in the following theorem.
\begin{thm}\label{roth}
Assume that there exists a strongly unbounded set of size $\non(\SMZ)$.
Then there exist a strongly null set of reals $X$ and a continuous
image $Y$ of $X$ such that $Y$ is not strongly null.
\end{thm}
\begin{proof}
Let $\kappa=\non(\SMZ)$.
Then there exist:
A strongly unbounded set $A=\{f_\alpha:\alpha<\kappa\}$,
and a set of reals $Y=\{y_\alpha:\alpha<\kappa\}$ that is not strongly null.
By standard translation arguments (see, e.g., \cite{prods}) we may assume
that $Y\sbst[0,1]$ and therefore think of $Y$ as a subset of $\Cantor$
($\Cantor$\index{$\Cantor$} is \Emph{Cantor's space}, which is endowed with the
product topology).
Consequently, the set $A'=\{f_\alpha' : \alpha<\kappa\}$,
where for each $\alpha$, $f_\alpha'(n) = 2f_\alpha(n)+y_\alpha(n)$ for all $n$,
is also strongly unbounded, and the mapping $A'\to Y$ defined
by $f(n) \mapsto f(n)\bmod 2$ is continuous and surjective.

It remains to show that $A'$ is a continuous image of a strongly null set of reals $X$.
Let $\Psi:\NN\to\R\sm\Q$ be a homeomorphism (e.g., taking continued fractions),
and let $X=\Psi[A']$. We claim that $X$ is strongly null.
Indeed, assume that $\seq{\epsilon_n}$ is a sequence of positive reals.
Enumerate $\Q=\{q_n : n\in\N\}$,
and choose for each $n$ an open interval $I_{2n}$ of length less than $\epsilon_{2n}$
such that $q_n\in I_{2n}$.
Let $U = \Union_n I_{2n}$.
Then $\R\sm U$ is $\sigma$-compact, thus $B=\Psi\inv[\R\sm U]=\NN\sm\Psi\inv[U]$
is a $\sigma$-compact and therefore bounded subset of $\NN$.
As $A'$ is strongly unbounded, $|A'\sm\Psi\inv[U]|=|A'\cap B|<\kappa=\non(\SMZ)$,
thus $\Psi[A'\cap B]=\Psi[A']\sm U$ is strongly null, so we can find intervals $I_{2n+1}$
of diameter at most $\epsilon_{2n+1}$
covering this set, and note that
$$X=\Psi[A']\sbst(\Psi[A']\sm U)\cup U\sbst\Union_{n\in\N}I_n.\qedhere$$
\end{proof}

From Theorem \ref{roth} it is possible to deduce that $\SMZ$ is not
provably closed under homeomorphic images.
The argument in the following proof, that is probably similar to
Rothberger's, was pointed out to us by T.\ Weiss.
Observe that $\SMZ$ is hereditary (that is, if $X$ is strongly null
and $Y$ is a subset of $X$, then $Y$ is strongly null too), and that it
is preserved under uniformly continuous images.

\begin{prop}\label{weiss}
If a hereditary property $P$ is not preserved under continuous images,
but is preserved under \emph{uniformly} continuous images,
then it is not preserved under homeomorphic images.
\end{prop}
\begin{proof}
Assume that $X$ satisfies $P$, $Y$ does not satisfy $P$,
and $f:X\to Y$ is a continuous surjection.
Let $\tilde X\sbst X$ (so that $\tilde X$ satisfies $P$)
be such that $f:\tilde X\to Y$ is a (continuous) bijection.
Then the set $f\sbst X\x Y$ (we identify $f$ with its graph)
is homeomorphic to $\tilde X$ which satisfies $P$, but the projection of $f$
on the second coordinate, which is a uniformly continuous image of $f$, is equal to $Y$.
Thus $f$ does not satisfy $P$.
\end{proof}

\subsection{Rothberger's property}\label{rothp}
This lead Rothberger to introduce the following topological version of strong measure zero
(which is preserved under continuous images).
Again, let $\scrA$ and $\scrB$ be collections of open covers of a topological space $X$.
Consider the following prototype of a selection hypothesis.
\bi
\item[$\sone(\scrA,\scrB)$\index{$\sone(\scrA,\scrB)$}:]
For each sequence $\seq{\cU_n}$ of members of $\scrA$,
there exist members $U_n\in\cU_n$, $n\in\N$, such that $\seq{U_n}\in\scrB$.
\ei

\begin{figure}[!htp]
\begin{changemargin}{-2cm}{-2cm}
\begin{center}
\myeps{\epsfysize=6 truecm {\epsfbox{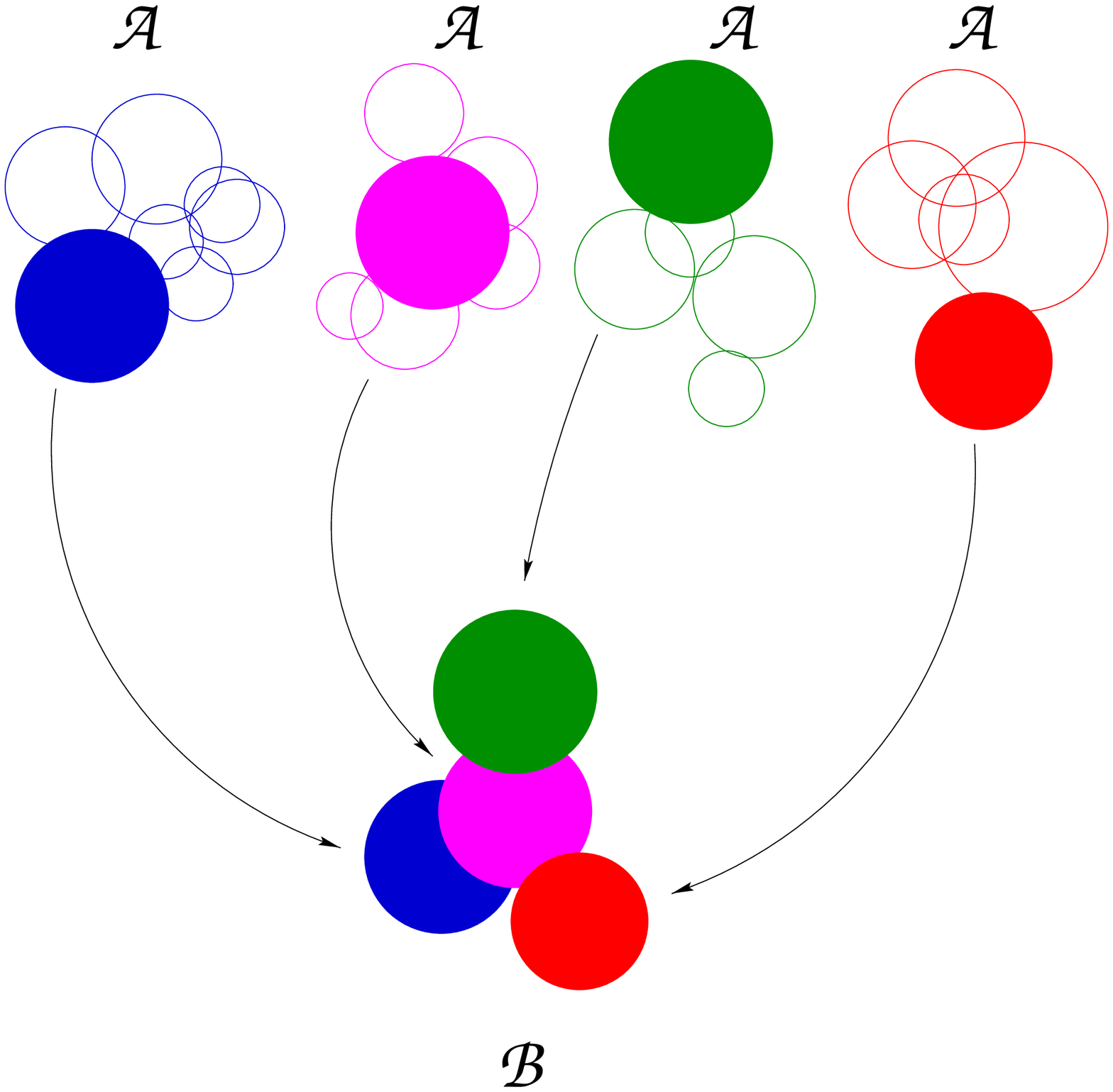}}}
\caption{$\sone(\scrA,\scrB)$}
\end{center}
\end{changemargin}
\end{figure}

Then Rothberger introduced the case $\scrA=\scrB=\O$ (the collection of all
open covers).\footnote{Originally, Rothberger denoted this property $C''$,
the reason being as follows.
In his 1919 paper \cite{Borel},
Borel considered several properties, which he enumerated as $A$, $B$, $C$, and so on.
The property that was numbered $C$ was that of strong measure zero.
\index{property $C$}.
Thus, Rothberger used
$C'$\index{property $C'$} to denote continuous images of elements of $C$, and
$C''$\index{property $C''$} to be what we now
call $\sone(\O,\O)$, since it implies $C'$.
}
Clearly, \emph{Rothberger's property}\index{Rothberger's property $\sone(\O,\O)$}
$\sone(\O,\O)$ implies being strongly null,
and the usual argument shows that every Luzin set $L$
satisfies $\sone(\O,\O)$.
Moreover, Fremlin and Miller \cite{FM} proved that for a metric space $\<X,d\>$,
$\sone(\O,\O)$ is the same as having strong measure zero with respect
to all metrics which generate the same topology as the one defined by $d$.

The question of the consistency of Borel's Conjecture was settled in
1976, when Laver in his deep work \cite{LAVER} showed that Borel's Conjecture is consistent.
We will return to Borel's Conjecture in Section \ref{BCrevisited}.

\section{Classification}

\subsection{More properties}
Having the terminology introduced thus far, we can also consider
the properties
$\sone(\Gamma,\O)$, $\sone(\Gamma,\Gamma)$, and $\sone(\O,\Gamma)$.
The last property turns out trivial (consider an open cover with no
$\gamma$-subcover),
but the first two make sense
even in the restricted setting of sets of reals.
These properties turn out much more restrictive than
Menger's property $\ufin(\O,\O)$, but they do not
admit an analogue of the Borel conjecture.
In fact, the set $H$\index{Hurewicz set $H$} from Theorem \ref{H} satisfies
$\sone(\Gamma,\O)$ (by an argument similar to that in
Theorem \ref{roth}) \cite{ideals}, and in fact
there always exist uncountable elements satisfying $\sone(\Gamma,\Gamma)$
\cite{coc2, wqn}.

\begin{prob}
~\be
\item (Bartoszy\'nski, et.\ al.\ \cite{ideals})
Does the set $H$ from Theorem \ref{H} satisfy $\sone(\Gamma,\Gamma)$?
\item Does there always exist a set of size $\b$ satisfying $\sone(\Gamma,\Gamma)$?
\ee
\end{prob}

\subsection{$\w$-covers}
We need not stop here, and may wish to consider other important types of covers which
appeared in the literature.
An open cover $\cU$ of $X$ is an
\Emph{$\omega$-cover} of $X$ if no single member of $\cU$ covers $X$,
but for each finite $F\sbst X$ there exists a single member of $\cU$ covering $F$.
Let $\Omega$\index{$\Omega$} denote the collection of open $\omega$-covers of $X$.
Then $\sone(\Omega,\Gamma)$ is equivalent to the \emph{$\gamma$-property}\index{Gerlits-Nagy $\gamma$-property}
introduced by Gerlits and Nagy (1982) in \cite{GN}, and
$\sone(\Omega,\Omega)$\index{Sakai property $\sone(\Omega,\Omega)$} was studied by Sakai (1988) in \cite{Sakai},
both properties naturally arising in the study of the space of continuous
real valued functions on $X$ (we will return to this in Section \ref{fspaces}).

\subsection{Arkhangel'ski\v{i}'s property}
Another prototype for a selection hypothesis generalizes a property studied by Arkhangel'ski\v{i},
also in the context of function spaces, in 1986 \cite{Arkhan86}.
The prototype is similar to $\ufin(\scrA,\scrB)$, but we do not ``glue''
the finite subcollections.
\bi
\item[$\sfin(\scrA,\scrB)$\index{$\sfin(\scrA,\scrB)$}:]
For each sequence $\seq{\cU_n}$
of members of $\scrA$, there exist finite (possibly empty)
subsets $\cF_n\sbst\cU_n$, $n\in\N$, such that $\Union_{n\in\N}\cF_n\in\scrB$.
\ei
Then the property studied by Arkhangel'ski\v{i} is equivalent to $\sfin(\Omega,\Omega)$.
\index{Arkhangel'ski\v{i} property $\sfin(\Omega,\Omega)$}

\subsection{The Scheepers Diagram}
\begin{figure}[!htp]
\begin{changemargin}{-4cm}{-3cm}
\begin{center}
{\scriptsize
$\xymatrix@R=6pt{
&
&
& \sr{\ufin(\Gamma,\Gamma)}{\b}\ar[r]
& \sr{\ufin(\Gamma,\Omega)}{\d}\ar[rr]
& & \sr{\ufin(\Gamma,\O)}{\d}
\\
&
&
& \sr{\sfin(\Gamma,\Omega)}{\d}\ar[ur]
\\
& \sr{\sone(\Gamma,\Gamma)}{\b}\ar[r]\ar[uurr]
& \sr{\sone(\Gamma,\Omega)}{\d}\ar[rr]\ar[ur]
& & \sr{\sone(\Gamma,\O)}{\d}\ar[uurr]
\\
&
&
& \sr{\sfin(\Omega,\Omega)}{\d}\ar'[u][uu]
\\
& \sr{\sone(\Omega,\Gamma)}{\p}\ar[r]\ar[uu]
& \sr{\sone(\Omega,\Omega)}{\cov(\M)}\ar[uu]\ar[rr]\ar[ur]
& & \sr{\sone(\O,\O)}{\cov(\M)}\ar[uu]
}$
}
\caption{The Scheepers Diagram}\label{SchDiag}
\end{center}
\end{changemargin}
\end{figure}
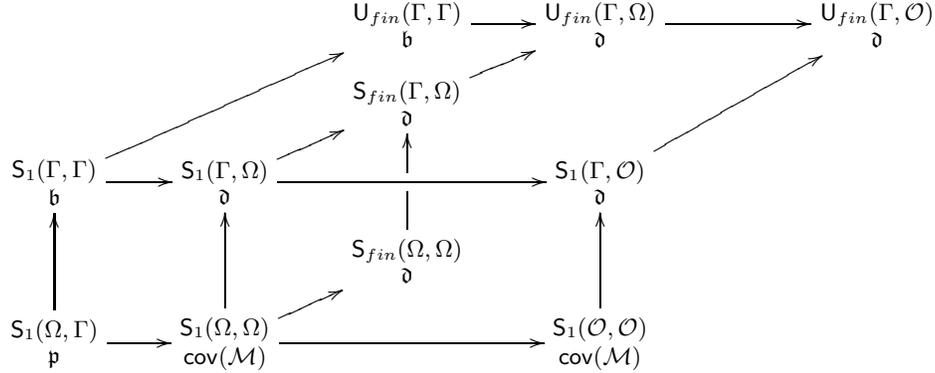

Thus far we have a selection hypothesis corresponding to each member of the $27$ element set
$\{\sone,\sfin,\ufin\}\x\{\Gamma,\Omega,\O\}^2$.
Fortunately, it suffices to consider only some of them.
First, observe that in the cases we consider,
$$\sone(\scrA,\scrB)\Impl\sfin(\scrA,\scrB)\Impl\ufin(\scrA,\scrB),$$
and we have the following
monotonicity property: For $\Pi\in\{\sone,\sfin,\ufin\}$, if $\scrA\sbst\scrC$ and $\scrB\sbst\scrD$, then:
$$\begin{array}{ccc}
\Pi(\scrA,\scrB) & \rightarrow & \Pi(\scrA,\scrD)\cr
\uparrow & & \uparrow\cr
\Pi(\scrC,\scrB) & \rightarrow & \Pi(\scrC,\scrD)
\end{array}$$
After removing trivial properties and proving equivalences
among the remaining ones (see \cite{coc2} for a summary of these),
we get the \Emph{Scheepers Diagram} (Figure \ref{SchDiag}).
In this diagram, as in the ones to follow, an arrow denotes implication,
and below each property we wrote its critical cardinality.
($\cov(\M)$\index{$\cov(\M)$} denotes the minimal cardinality of a cover of
$\R$ by meager sets, and $\p$ is the \emph{\psin{} number} to be defined in
Section \ref{MTP}. See \cite{BlassHBK} for information on these cardinals as
well as other cardinals which we mention later.)

There remain only two problems concerning this diagram.

\begin{prob}[{Just, Miller, Scheepers, Szeptycki \cite{coc2}}]\label{classify}
~\be
\item Does $\ufin(\Gamma,\Omega)$ imply $\sfin(\Gamma,\Omega)$?
\item If not, does $\ufin(\Gamma,\Gamma)$ imply $\sfin(\Gamma,\Omega)$?
\ee
\end{prob}

All other implications are settled in \cite{coc1, coc2},
using two methods. One approach uses
consistency results concerning the values of the critical cardinalities.
For example, it is consistent that $\b<\d$, thus none of the properties with
critical cardinality $\d$ can imply any of those with critical cardinality $\b$.
Another approach is by transfinite constructions under \CH{}
(such as special kinds of Luzin and Sierpi\'nski sets).
Recently, an approach combining these two approaches was investigated --
see, e.g., \cite{JORG, BartHBK, gamma7}.

\subsection{Borel covers}
\begin{figure}[!htp]
\begin{changemargin}{-4cm}{-3cm}
\begin{center}
{\tiny
$\xymatrix@C=-2pt@R=6pt{
&
&
& \sr{\ufin(\Gamma,\Gamma)}{\b}\ar[rr]\ar@{.>}[dr]^?
&
& \sr{\ufin(\Gamma,\Omega)}{\d}\ar[rrrrr]\ar@/_/@{.>}[dl]_?
&
&
&
&
&
&
& \sr{\ufin(\Gamma,\O)}{\d}
\\
&
&
&
& \sr{\sfin(\Gamma,\Omega)}{\d}\ar[ur]
\\
& \sr{\sone(\Gamma,\Gamma)}{\b}\ar[rr]\ar[uurr]
&
& \sr{\sone(\Gamma,\Omega)}{\d}\ar[rrr]\ar[ur]
&
&
& \sr{\sone(\Gamma,\O)}{\d}\ar[uurrrrrr]
\\
  \sr{\sone(\BG,\BG)}{\b}\ar[ur]\ar[rr]
&
& \sr{\sone(\BG,\BO)}{\d}\ar[ur]\ar[rrr]
&
&
& \sr{\sone(\BG,\B)}{\d}\ar[ur]
\\
&
&
&
& \sr{\sfin(\Omega,\Omega)}{\d}\ar'[u]'[uu][uuu]
\\
\\
&
& \sr{\sfin(\BO,\BO)}{\d}\ar[uuu]\ar[uurr]
\\
& \sr{\sone(\Omega,\Gamma)}{\p}\ar'[r][rr]\ar'[uuuu][uuuuu]
&
& \sr{\sone(\Omega,\Omega)}{\cov(\M)}\ar'[uuuu][uuuuu]\ar'[rr][rrr]\ar[uuur]
&
&
& \sr{\sone(\O,\O)}{\cov(\M)}\ar[uuuuu]
\\
  \sr{\sone(\BO,\BG)}{\p}\ar[uuuuu]\ar[rr]\ar[ur]
&
& \sr{\sone(\BO,\BO)}{\cov(\M)}\ar[uu]\ar[ur]\ar[rrr]
&
&
& \sr{\sone(\B,\B)}{\cov(\M)}\ar[uuuuu]\ar[ur]
}$
}
\caption{The extended Scheepers Diagram}\label{extSch}
\end{center}
\end{changemargin}
\end{figure}
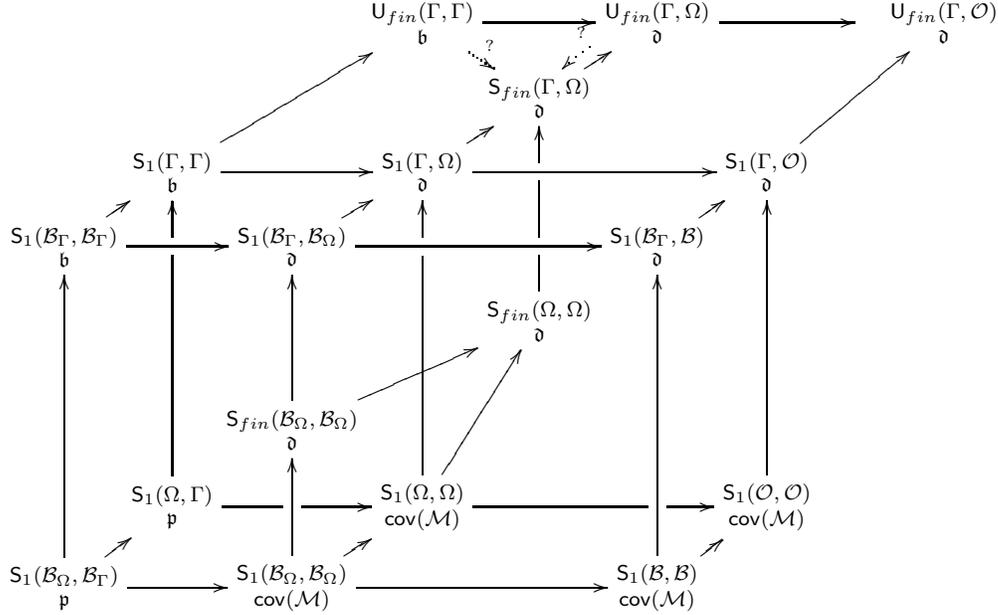

There are other natural types of covers which appear in the literature,
but probably the first natural question is: What happens if we replace
``open'' by ``countable Borel'' in the types of covers which we consider?
Let $\B,\BO,\BG$\index{$\B$}\index{$\BO$}\index{$\BG$} denote the collections of
\emph{countable Borel} covers,
$\w$-covers, and $\gamma$-covers of the given space, respectively.
It turns out that the same analysis is applicable
when we plug in these families instead of $\O,\Omega,\Gamma$,
and in fact one gets more equivalences. Moreover, some of the resulting properties turn out
equivalent to properties which appeared in the literature in other guises.
This is shown in \cite{CBC}, where it is also shown that no arrows can be added
(except perhaps those corresponding to Problem \ref{classify})
to the extended diagram (Figure \ref{extSch}).

\subsection{Borel's Conjecture revisited}\label{BCrevisited}
It is easy to verify that every countable set of reals $X$ satisfies
the strongest property in the extended Scheepers Diagram \ref{extSch},
namely, $\sone(\BO,\BG)$. By Laver's result mentioned in Section \ref{rothp},
we know that it is consistent that all properties between $\sone(\BO,\BG)$
and $\sone(\O,\O)$ (inclusive) are consistent to hold only for countable sets
of reals.

All other classes in the original Scheepers Diagram \ref{SchDiag}
contain uncountable elements:
Recall that every $\sigma$-compact set satisfies the Hure\-wicz property $\ufin(\O,\Gamma)$.
Now, $\sfin(\Omega,\Omega)$ is equivalent to satisfying Menger's property $\ufin(\O,\O)$
in all finite powers \cite{coc2}. Thus, every $\sigma$-compact set satisfies $\sfin(\Omega,\Omega)$.
As for the remaining properties, recall from Section \ref{ZFCexamples} that $\sone(\Gamma,\Gamma)$
always contains an uncountable element.
In other words, none of the properties except those mentioned in the
previous paragraph can satisfy an analogue of
Borel's Conjecture. However, by a result of Miller,
Borel's Conjecture for $\sone(\BG,\BG)$ is consistent \cite{ideals}.

\begin{prob}\label{domiBC}
~\be
\item (folklore) Is it consistent that every set of reals which satisfies
$\sone(\BG,\B)$ is countable?
\item What about $\sfin(\BO,\BO)$ and $\sone(\BG,\BO)$?
\ee
\end{prob}
A combinatorial formulation of the first question in Problem \ref{domiBC}
is obtained by replacing $\sone(\BG,\B)$ with the equivalent property
``every Borel image in $\NN$ is not dominating''.

Other interesting investigations in this direction are
of the form: Is it consistent that a certain property in the diagram
satisfies Borel's Conjecture, whereas another one does not?
Some results in this direction are the following.
\begin{thm}[Miller \cite{MillerBC}]\label{miller}
~\be
\item Borel's Conjecture for $\sone(\O,\O)$ implies
Borel's Conjecture (for strong measure zero);
\item Borel's Conjecture for $\sone(\Omega,\Gamma)$ does not imply Borel's Conjecture.
\ee
\end{thm}

The proofs use, of course, combinatorial arguments (and forcing in the second
case: The model is obtained by adding $\aleph_2$ dominating reals with a finite support iteration to
a model of \CH{}).
Using yet more combinatorial arguments, it is possible to extend Theorem \ref{miller}.

\begin{thm}[Weiss, et.\ al.\ \cite{prods}]
Borel's Conjecture for $\sone(\Omega,\Omega)$ implies Borel's Conjecture.
Consequently, Borel's Conjecture for $\sone(\Omega,\Gamma)$ does not imply Borel's Conjecture
for $\sone(\Omega,\Omega)$.
\end{thm}

This settles completely this investigation when we restrict attention to
the original Scheepers Diagram \ref{SchDiag}.

\section{Preservation of properties}

\subsection{Continuous images}
It is easy to see that all properties in the Scheepers Diagram \ref{SchDiag}
(as well as all other selection properties in this paper)
are preserved under continuous images \cite{coc2}.
Similarly, the selection properties involving Borel covers are preserved under Borel
images \cite{CBC}. The situation is not as good concerning other types of preservation\dots

\subsection{Additivity}
In \cite{coc2} (1996), Just, Miller, Scheepers, and Szeptycki raised the following
\Emph{additivity problem}: It is easy to see that some of the properties
in the Scheepers diagram are (provably) preserved under taking finite and even countable
unions (i.e., they are \Emph{$\sigma$-additive}). What about the remaining ones?
Figure \ref{add}(a) summarizes the knowledge that was available at the point the
question was posed, where the positions are according to Figure \ref{SchDiag}.

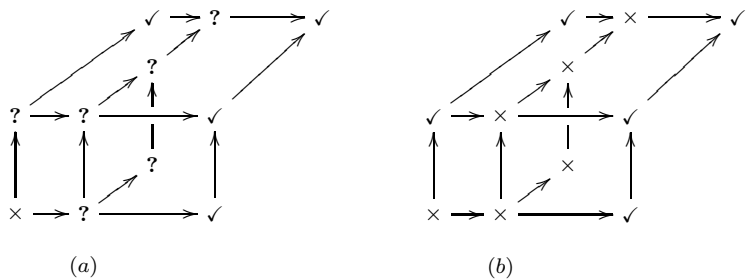
\begin{figure}[!htp]
\begin{changemargin}{-4cm}{-3cm}
\begin{center}
{\tiny
$\xymatrix@R=6pt@C=10pt{
&
&
& \txt{\checkmark}\ar[r]
& \txt{\textbf{?}}\ar[rr]
& & \txt{\checkmark}
\\
&
&
& \txt{\textbf{?}}\ar[ur]
\\
& \txt{\textbf{?}}\ar[r]\ar[uurr]
& \txt{\textbf{?}}\ar[rr]\ar[ur]
& & \txt{\checkmark}\ar[uurr]
\\
&
&
& \txt{\textbf{?}}\ar'[u][uu]
\\
& \times\ar[r]\ar[uu]
& \txt{\textbf{?}}\ar[uu]\ar[rr]\ar[ur]
& & \txt{\checkmark}\ar[uu]\\
&                    & (a)
}$
\quad
$\xymatrix@R=6pt@C=10pt{
&
&
& \txt{\checkmark}\ar[r]
& \times\ar[rr]
& & \txt{\checkmark}
\\
&
&
& \times\ar[ur]
\\
& \txt{\checkmark}\ar[r]\ar[uurr]
& \times\ar[rr]\ar[ur]
& & \txt{\checkmark}\ar[uurr]
\\
&
&
& \times\ar'[u][uu]
\\
& \times\ar[r]\ar[uu]
& \times\ar[uu]\ar[rr]\ar[ur]
& & \txt{\checkmark}\ar[uu]\\
&                    & (b)
}$
}
\end{center}
\caption[The additivity problem and its solution]{The additivity problem (a) and its solution (b)}
\end{changemargin}
\end{figure}

In 1999, Scheepers \cite{wqn}
proved that the answer is positive for $\sone(\Gamma,\Gamma)$.

The key observation towards solving the problem for the remaining properties
was the following analogue of Hurewicz' Theorem \ref{hure}.
According to Blass, a family $Y\sbst\NN$ is \Emph{finitely dominating} if
for each $g\in\NN$ there exist $k$ and $f_1,\dots,f_k\in Y$
such that $g(n)\le\max\{f_1(n),\dots,f_k(n)\}$ for all
but finitely many $n$.

\begin{thm}[\cite{huremen1}]\label{key}
A set of reals $X$ satisfies $\ufin(\Gamma,\Omega)$
if, and only if, each continuous
image of $X$ in $\NN$ is not finitely dominating.
\end{thm}

Motivated by this observation, the following theorem was proved.
(As usual, $\c=2^{\aleph_0}$\index{$\c$} denotes the size of the continuum.)
\begin{thm}[Bartoszy\'nski, Shelah, et.\ al.\ \cite{huremen2, AddQuad}]\label{add}
Assume \CH{} (or just $\cov(\M)=\c$). Then
there exist sets of reals $L_0$, and $L_1$ satisfying
$\sone(\BO,\BO)$ such that
$L_0\cup L_1$ is finitely dominating.
\end{thm}
The proof used a tricky ``power sharing'' between $L_0$ and $L_1$ during
their transfinite-inductive construction.

Consequently, \emph{none of the remaining properties} is provably preserved
under taking finite unions (Figure \ref{add}(b)). This also settled all corresponding problems
in the Borel case (see the extended diagram \ref{extSch}).
Interestingly, the simple observation in Theorem \ref{key} was also
the key behind proving that consistently (namely, assuming NCF),
$\ufin(\O,\Omega)$ is $\sigma$-additive.

\subsubsection*{Bad transmission of knowledge}
The transmission of knowledge concerning the additivity problem was very poor
(take a deep breathe):
It posteriorly turns out that
if we restrict attention to the open case only, then
the additivity problem was already implicitly solved earlier.
In 1999, Scheepers \cite{lengthdiags} constructed
sets of reals $L_0$, and $L_1$ satisfying
$\sone(\Omega,\Omega)$ such that
$L_0+L_1$ is finitely dominating. It is easy to see that this implies
that $L_0\cup L_1$ is finitely dominating, which settles the problem
if we add the missing ingredient
Theorem \ref{key}, which, funnily, seems to be the simplest part of the solution.
Scheepers was unaware of this observation and consequently of that solving the
additivity problem completely, but he did point out that his construction implied
that the properties between $\sone(\Omega,\Omega)$ and $\sfin(\Omega,\Omega)$
(inclusive) are not provably additive. In turn, we were not aware of this
when we proved Theorem \ref{add}.
On top of that, Theorem \ref{key}, the open part of Theorem \ref{add},
and the result concerning NCF were independently proved
by Eisworth and Just in \cite{gamma7},
and similar results were also independently obtained by Banakh, Nickolas, and Sanchis
in \cite{BNS}. As if this is not enough, the main ingredient of the result
concerning NCF was also independently obtained by Blass \cite{BlassNA}.

This is a good point to recommend the reader announce his new results in the
\SPMBul{} (see \cite{spm}), so as to avoid similar situations.

\subsection{Hereditarity}

A property (or a class of topological spaces) is
\Emph{hereditary} if it is preserved under taking subsets.
Despite the fact that the properties in the Scheepers diagram
are (intuitively) notions of smallness,
none of them is (provably) hereditary.
Define a topology on the space $P(\N)$\index{$P(\N)$}
of all sets of natural numbers by identifying it
with $\Cantor$. Note that $\fin$\index{$\fin$},
the collection of finite subsets of $\N$, is dense in $P(\N)$.
By taking increasing enumerations, we have that the
\emph{Rothberger space}\index{Rothberger space $\roth$}
$\roth=P(\N)\sm\fin$ is identified with the space of increasing sequences of natural numbers,
which in turn is homeomorphic to the Baire space $\NN$.

\begin{thm}[{Babinkostova, Ko\v{c}inac, and Scheepers \cite{coc8}; Bartoszy\'nski, et.\ al.\ \cite{ideals}}]
Assuming (a portion of) \CH,
there exists a set of reals $X$ satisfying $\sone(\Omega,\Gamma)$
and a countable subset $Q$ of $X$,
such that $X\sm Q$ does not satisfy $\ufin(\O,\O)$.
\end{thm}
The proof for this assertion is a modification of the construction
of Galvin and Miller \cite{GM}, with $X\sbst P(\N)$, $Q=\fin$, so that
$X\sm\fin\sbst\roth$ is dominating (when viewed as a subset of
$\NN$) which is what we look for in light of Hurewicz' Theorem \ref{hure}.

However, in the Borel case,
$\sone(\B,\B)$ and all properties of the
form $\Pi(\BG,\scrB)$ are hereditary.
This is immediate from the combinatorial characterizations \cite{CBC},
but is also easy to prove directly  \cite{ideals}.
However, not all properties in the Borel case are hereditary, e.g.,
Miller \cite{MilNonGamma} proved that assuming \CH{},
$\sone(\BO,\BG)$ is not hereditary.

\begin{prob}[{Bartoszy\'nski, et.\ al.\ \cite{ideals}}]
Is any of $\sone(\BO,\BO)$ or $\sfin(\BO,\BO)$ hereditary?
\end{prob}

\section{The Minimal Tower Problem and $\tau$-covers}

\subsection{``Rich'' covers of spaces}\footnote{Richness can be viewed as
some sort of redundancy, and some people are richer than others.
This is precisely the case with
\Emph{rich covers}.}
To avoid trivialities, let us decide that
when we say that $\cU$ is a \Emph{cover} of $X$,
we mean that $X=\cup\cU$ (the usual requirement), and
in addition no single member of $\cU$ covers $X$.
We also always assume that the space $X$ is infinite.

$\cU$ is a \Emph{large cover} of $X$ if for each $x\in X$, there exist infinitely many
$U\in\cU$ such that $x\in U$.
Let $\Lambda$\index{$\Lambda$} denote the collection of large open covers of $X$.
Then
$$\Gamma\sbst\Omega\sbst\Lambda,$$
the last implication being a cute exercise.

$\cU$ is a \Emph{$\tau$-cover} of $X$ if it is a large cover of $X$,
and for each $x,y\in X$, at least one of the sets
$\{U\in\cU : x\in U\mbox{ and }y\nin U\}$ or
$\{U\in\cU : y\in U\mbox{ and }x\nin U\}$ is finite.
This is a reminiscent of the negation of the $T_1$ property,
where here the notion is applied ``modulo finite''.
Let $\Tau$\index{$\Tau$} (uppercase $\tau$)
denote the collection of open $\tau$-covers of $X$.
Then
$$\Gamma\sbst\Tau\sbst\Omega.$$
The motivation behind the definition of $\tau$-covers is the
Minimal Tower Problem.

\subsection{The Minimal Tower Problem}\label{MTP}
Recall that $\roth$ is the collection of infinite subsets of $\N$.
$A\as B$ means that $A\sm B$ is finite.
$\cF\sbst\roth$ is \Emph{centered} if
for each finite $\cA\sbst\cF$, $\cap\cA$ is infinite.
$A\in\roth$ is a \Emph{\psin{}} of
$\cF$ if $A\as B$ for all $B\in\cF$.
Let $\p$\index{$\p$} denote the minimal size of a centered family $\cF\sbst\roth$ with no \psin{},
and $\t$\index{$\t$} denote the minimal size of a $\as$-linearly ordered family $\cF\sbst\roth$
which has no \psin{}.
Then $\p\le\t$, and the \Emph{Minimal Tower Problem} is:

\begin{prob}
Is it provable that $\p=\t$?
\end{prob}

This is one of the major and oldest problems of infinitary
combinatorics. Allusions to this problem can be
found in Rothberger's 1940's works (see, e.g., \cite{ROTH48}).
The problem is explicitly mentioned in van Douwen's survey
\cite{vD}, and quoted in Vaughan's survey \cite[Problem 333]{Vaughan},
where it is considered the most interesting open problem in the
field. Extensive work of Shelah and others in the field solved
all problems of this sort which were mentioned in \cite{vD},
except for the Minimal Tower Problem.

\subsection{A dictionary}\label{dictionary}
How is this problem related to $\tau$-covers?
Each type of rich cover corresponds to some combinatorial notion
(of richness), as follows.
If we confine attention to sets of reals, then
we may assume that all open covers we consider are countable.
Thus, assume that $\cU=\seq{U_n}$ is a countable (not necessarily open)
cover of $X$, and define the \Emph{Marczewski characteristic function} of $\cU$
\cite{marczewski38} by
$$h(x) = \{ n : x\in U_n\},$$
so that $h:X\to P(\N)$.
$h$ is continuous if the sets $U_n$ are clopen,
and Borel if the sets $U_n$ are Borel.

We have the following \Emph{dictionary} translating properties
of $\cU$ to properties of the image $h[X]$.
\begin{center}
\begin{tabular}{|c|c|}
\hline
$\seq{U_n}$    & $h[X]$\\
\hline\hline
large cover    & subset of $\roth$\\
\hline
$\w$-cover     & centered\\
\hline
$\tau$-cover   & linearly ordered by $\as$\\
\hline
$\gamma$-cover & cofinite sets\\
\hline
contains a $\gamma$-cover & has a pseudo-intersection\\
\hline
\end{tabular}
\end{center}
(Note that, since we assume that $X\not\sbst U_n$ for all $n$, we have that
$h[X]$ is centered if, and only if, it is a base for a nonprincipal filter on $\N$.)

\subsection{Topological approximations of the Minimal Tower Problem}
For families $\scrB\sbst\scrA$ of covers of a space $X$, define the property
\Emph{$\scrA$ choose $\scrB$} as follows.
\bi
\item[$\binom{\scrA}{\scrB}$\index{$\binom{\scrA}{\scrB}$}:] For each $\cU\in\scrA$ we can choose a subset $\cV\sbst\cU$ such that $\cV\in\scrB$.
\ei
This is a prototype for many classical topological notions, most notably compactness and being Lindel\"of.

In 1982, Gerlits and Nagy \cite{GN} introduced the
\Emph{$\gamma$-property} $\binom{\Omega}{\Gamma}$, and proved
that it is equivalent to $\sone(\Omega,\Gamma)$.

By the dictionary (Section \ref{dictionary}), we have that
\be
\item $\non\binom{\Omega}{\Gamma}=\non\binom{\BO}{\BG}=\p$; and
\item $\non\binom{\Tau}{\Gamma}=\non\binom{\BT}{\BG}=\t$.
\ee

Clearly, $\binom{\Omega}{\Gamma}\sbst\binom{\Tau}{\Gamma}$, and
$\binom{\BO}{\BG}\sbst\binom{\BT}{\BG}$. We therefore have the following
topological problems related to the Minimal Tower Problem.
\be
\item Is $\binom{\Omega}{\Gamma}=\binom{\Tau}{\Gamma}$ ?
\item Is $\binom{\BO}{\BG}=\binom{\BT}{\BG}$ ?
\ee
Observe that the answer is ``No'' for both problems if $\p<\t$ is consistent.

But it turns out that both problems can be solved outright in ZFC:
The first problem was solved by Shelah \cite{tau}, who proved
that $\Cantor$ satisfies $\binom{\Tau}{\Gamma}$.
A set satisfying $\binom{\Omega}{\Gamma}$
must satisfy $\sone(\O,\O)$ and therefore have strong measure zero.
Obviously, $\Cantor$ does not have strong measure zero (it has $[0,1]$ as
a uniformly continuous image).
A modification of Shelah's construction yields
a negative solution to the second question as well: Indeed,
$\NN$ satisfies $\binom{\Tau}{\Gamma}$ \cite{tautau}, and
Borel images of $\NN$ are analytic, and can therefore
be represented as continuous images of $\NN$, so Borel images
of $\NN$ satisfy $\binom{\Tau}{\Gamma}$, and therefore
$\NN$ satisfies $\binom{\BT}{\BG}$ \cite{tautau}.

\subsection{Tougher topological approximations}
It is easy to see that $\non(\sone(\Omega,\Gamma))=\non(\sone(\BO,\BG))=\p$,
and $\non(\sone(\Tau,\Gamma))=\non(\sone(\BT,\allowbreak\BG))=\t$ \cite{tau},
and we have the following implications:
$$\begin{array}{ccc}
\sone(\Omega,\Gamma) & \rightarrow & \sone(\Tau,\Gamma)\cr
\uparrow &                         & \uparrow\cr
\sone(\BO,\BG) & \rightarrow & \sone(\BT,\BG)
\end{array}$$
We therefore arrive at the following ``tighter'' approximations.
\be
\item Is $\sone(\Omega,\Gamma)=\sone(\Tau,\Gamma)$ ?
\item Is $\sone(\BO,\BG)=\sone(\BT,\BG)$ ?
\ee
Again, $\p<\t$ implies a negative answer for both problems.
It turns out that the answer is indeed negative,
at least assuming \CH{} \cite{tautau}.
The main ingredients in the proof of this are the following:
$\sone(\BT,\BG)=\sone(\BG,\BG)\cap\binom{\BT}{\BG}$, and is therefore countably additive.
Now, \CH{} implies that there exists a hereditary $\sone(\BO,\BG)$-set $X\sbst [0,1]$
(e.g., Brendle \cite{JORG} or Miller \cite{MilNonGamma}).
By a result of Galvin and Miller (1984) \cite{GM},
if $Y\sbst X$ is not $F_\sigma$ or not $G_\delta$, then
$(X\sm Y)\cup (Y+1)\nin\sone(\Omega,\Gamma)$.
This solves both problems at once.

The striking observation we get from this journey of topological
approximations to the Minimal Tower Problem is that, despite
the fact that the purely combinatorial problem is very difficult,
if we add to it topological structure we get a negative answer
quite easily.

It is surprising, though, that the following related problem
remains open.
\begin{prob}[\cite{tautau}]
Is $\binom{\Omega}{\Gamma}=\binom{\Omega}{\Tau}$?
\end{prob}
The critical cardinality of $\binom{\Omega}{\Tau}$ is $\p$ \cite{ShTb768},
so both properties have the same critical cardinality.

\begin{rem}
The reader interested in another topological study which was
inspired by (and is related to) the Minimal Tower Problem is
referred to Machura's \cite{funcpt}.
\end{rem}

\subsection{Known implications and critical cardinalities}

Having this new notion of rich covers, namely $\tau$-covers,
it is interesting to try and add it to Scheepers' framework
of selection principles.
This yields many more properties, but again some of them can
be proved to be equivalent \cite{tautau}.

As already mentioned,
a basic tool to prove \emph{nonimplications}\index{nonimplications}
is the computation of
critical cardinalities: If $P$ and $Q$ are properties with
$\non(P)<\non(Q)$ consistent, then $Q$ does not imply $P$.
The critical cardinalities of most of the new properties were found
in \cite{tautau} and in Shelah, et.\ al.\ \cite{ShTb768}.
Still, $6$ critical cardinalities remained unsettled.
These remaining cardinalities were addressed by
Mildenberger, Shelah, et.\ al., \cite{MShT:858}.
We give here one example of that treatment, and quote the main results.
(Everything in the remainder of this subsection is quoted, without further notice,
from \cite{MShT:858}).

Let $\CG$, $\CT$, and $\CO$ denote the collections of
\emph{clopen} $\gamma$-covers, $\tau$-covers, and $\w$-covers of $X$, respectively.
Recall that, since we are dealing with sets of reals, we may assume that all open covers
are countable. Restricting attention to countable covers,
we have the following, where an arrow denotes inclusion:
$$\begin{matrix}
\BG      & \impl & \BT      & \impl & \BO      & \impl & \B      \\
\uparrow &       & \uparrow &       & \uparrow &       & \uparrow \\
\Gamma   & \impl & \Tau     & \impl & \Omega   & \impl & \O  \\
\uparrow &       & \uparrow &       & \uparrow &       & \uparrow \\
\CG      & \impl & \CT      & \impl & \CO      & \impl & C
\end{matrix}$$
As each of the properties $\Pi(\cdot\ ,\cdot)$, $\Pi\in\{\sone,\sfin,\ufin\}$,
is monotonic in its first variable, we have that for each $x,y\in\{\Gamma, \Tau, \Omega,\O\}$,
$$\Pi(\B_x,\B_y)\impl\Pi(x,y)\impl\Pi(C_x,C_y)$$
(here $C_\O:=C$ and $\B_\O:=\B$).
Consequently,
$$\non(\Pi(\B_x,\B_y))\le\non(\Pi(x,y))\le\non(\Pi(C_x,C_y)).$$

\begin{defn}\label{diagble}
We use the short notation $\forall^\oo$ for ``for all but finitely many''
and $\exists^\oo$ for ``there exist infinitely many''.
\be
\item $A\in \arrays$ is a \emph{$\gamma$-array}
if $(\forall n)(\forall^\oo m)\ A(n,m)=1$.
\item $\cA\sbst\arrays$ is a \emph{$\gamma$-family} if each $A\in\cA$ is a
$\gamma$-array.
\item A family $\cA\sbst\arrays$ is \emph{finitely $\tau$-diagonalizable} if
there exist finite (possibly empty) subsets $F_n\sbst\N$, $n\in\N$, such
that:
\be
\item For each $A\in\cA$: $(\exists^\oo n)(\exists m\in F_n)\ A(n,m)=1$;
\item For each $A,B\in\cA$:\\
\begin{tabular}{ll}
Either  & $(\forall^\oo n)(\forall m\in F_n)\ A(n,m)\le B(n,m)$,\\
or      & $(\forall^\oo n)(\forall m\in F_n)\ B(n,m)\le A(n,m)$.
\end{tabular}
\ee
\ee
\end{defn}

Using the dictionary of Section \ref{dictionary}, one can prove the following.
(Note that $\arrays$ is topologically the same as the Cantor space $\Cantor$.)

\begin{thm}\label{charSfinGT}
For a set of reals $X$, the following are equivalent:
\be
\item $X$ satisfies $\sfin(\BG,\BT)$; and
\item For each Borel function $\Psi:X\to\arrays$, if $\Psi[X]$ is a $\gamma$-family,
then it is finitely $\tau$-diagonalizable.
\ee
The corresponding assertion for $\sfin(\CG,\CT)$
holds when ``Borel'' is replaced by ``continuous''.
\end{thm}

\begin{lem}\label{dgbl=b}
The minimal cardinality of a $\gamma$-family which is not finitely $\tau$-diagonaliz\-able is $\b$.
\end{lem}
(To appreciate the result in Lemma \ref{dgbl=b}, the reader may wish to try proving it for a while.)
Having Theorem \ref{charSfinGT} and Lemma \ref{dgbl=b}, we get that
$$\b=\non(\sfin(\BG,\BT))\le\non(\sfin(\Gamma,\Tau))\le\non(\sfin(\CG,\CT))=\b,$$
and therefore $\non(\sfin(\Gamma,\Tau))=\b$.

Note that Theorem \ref{charSfinGT} reduced the original topological question
into a purely combinatorial question. Its solution in Lemma \ref{dgbl=b} may be
of independent interest to those working on this sort of pure combinatorics.

It follows that $\non(\sone(\Gamma,\Tau))=\b$, and using a similar approach,
it can be proved that $\non(\sone(\Tau,\Tau))=\t$, and $\non(\sfin(\Tau,\Tau))=\min\{\b,\s\}$.

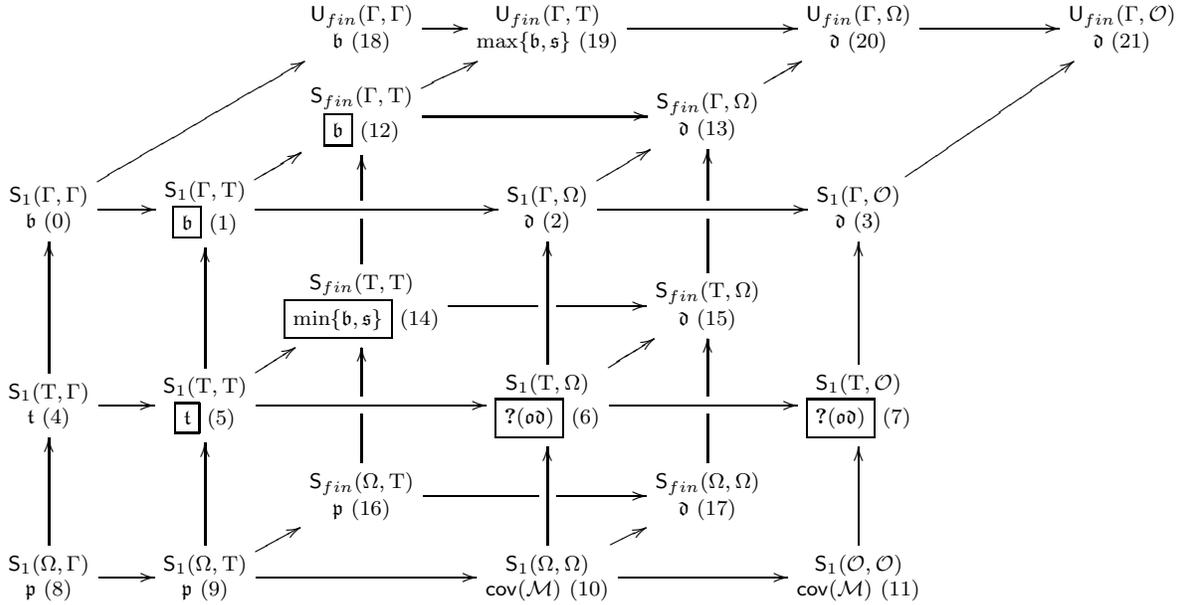
\begin{figure}[!ht]
\renewcommand{\sr}[2]{{\txt{$#1$\\$#2$}}}
{\tiny
\begin{changemargin}{-3cm}{-3cm}
\begin{center}
$\xymatrix@C=7pt@R=6pt{
&
&
& \sr{\ufin(\Gamma,\Gamma)}{\b~~ (18)}\ar[r]
& \sr{\ufin(\Gamma,\Tau)}{\max\{\b,\s\}~~ (19)}\ar[rr]
&
& \sr{\ufin(\Gamma,\Omega)}{\d~~ (20)}\ar[rrrr]
&
&
&
& \sr{\ufin(\Gamma,\O)}{\d~~ (21)}
\\
&
&
& \sr{\sfin(\Gamma,\Tau)}{\fb{\b}~~ (12)}\ar[rr]\ar[ur]
&
& \sr{\sfin(\Gamma,\Omega)}{\d~~ (13)}\ar[ur]
\\
\sr{\sone(\Gamma,\Gamma)}{\b~~ (0)}\ar[uurrr]\ar[rr]
&
& \sr{\sone(\Gamma,\Tau)}{\fb{\b}~~ (1)}\ar[ur]\ar[rr]
&
& \sr{\sone(\Gamma,\Omega)}{\d~~ (2)}\ar[ur]\ar[rr]
&
& \sr{\sone(\Gamma,\O)}{\d~~ (3)}\ar[uurrrr]
\\
&
&
& \sr{\sfin(\Tau,\Tau)}{\fb{\min\{\b,\s\}}~~ (14)}\ar'[r][rr]\ar'[u][uu]
&
& \sr{\sfin(\Tau,\Omega)}{\d~~ (15)}\ar'[u][uu]
\\
\sr{\sone(\Tau,\Gamma)}{\t~~ (4)}\ar[rr]\ar[uu]
&
& \sr{\sone(\Tau,\Tau)}{\fb{\t}~~ (5)}\ar[uu]\ar[ur]\ar[rr]
&
& \sr{\sone(\Tau,\Omega)}{\fbox{\textbf{?}$(\fo)$}~~ (6)}\ar[uu]\ar[ur]\ar[rr]
&
& \sr{\sone(\Tau,\O)}{\fbox{\textbf{?}$(\fo)$}~~ (7)}\ar[uu]
\\
&
&
& \sr{\sfin(\Omega,\Tau)}{\p~~ (16)}\ar'[u][uu]\ar'[r][rr]
&
& \sr{\sfin(\Omega,\Omega)}{\d~~ (17)}\ar'[u][uu]
\\
\sr{\sone(\Omega,\Gamma)}{\p~~ (8)}\ar[uu]\ar[rr]
&
& \sr{\sone(\Omega,\Tau)}{\p~~ (9)}\ar[uu]\ar[ur]\ar[rr]
&
& \sr{\sone(\Omega,\Omega)}{\cov(\M)~~ (10)}\ar[uu]\ar[ur]\ar[rr]
&
& \sr{\sone(\O,\O)}{\cov(\M)~~ (11)}\ar[uu]
}$
\end{center}
\end{changemargin}
}

\caption{The Scheepers diagram, enhanced with $\tau$-covers}\label{tauSch}
\end{figure}

What about the remaining two cardinals?
It is not difficult to see that $\non(\sfin(\Tau,\Omega))=\non(\sfin(\Tau,\O))$,
call this joint cardinal $\fo$, the \emph{$o$-diagonalization number}\index{$o$-diagonalization number};
the reason for this to be explained soon.
The surviving properties (in the open case)
appear in Figure \ref{tauSch}, with their critical cardinalities,
and serial numbers (for later reference).
The new cardinalities computed in \cite{MShT:858} are framed.

By Figure \ref{tauSch},
$$\cov(\M)=\non(\sone(\O,\O))\le\non(\sone(\Tau,\O))\le\non(\sone(\Gamma,\O))=\d,$$
thus $\cov(\M)\le\fo\le\d$.

\begin{defn}\label{odiagbl}
A family $\cA\sbst\arrays$ is a \emph{$\tau$-family}
if:
\be
\item For each $A\in\cA$: $(\forall n)(\exists^\oo m)\ A(n,m)=1$;
\item For each $A,B\in\cA$ and each $n$:\\
\begin{tabular}{ll}
Either & $(\forall^\oo m)\ A(n,m)\le B(n,m)$,\\
or     & $(\forall^\oo m)\ B(n,m)\le A(n,m)$.
\end{tabular}
\ee
A $\tau$-family $\cA$ is \emph{$o$-diagonalizable} if
there exists a function $g:\N\to\N$, such that:
$$(\forall A\in\cA)(\exists n)\ A(n,g(n))=1.$$
(Equivalently, $(\forall A\in\cA)(\exists^\oo n)\ A(n,g(n))=1$).
\end{defn}

As in Theorem \ref{charSfinGT}, $\sone(\BT,\B)$ and $\sone(\CT,C)$
have a natural combinatorial characterization.
This characterization implies that
$\fo$ is equal to the minimal cardinality of a $\tau$-family that is not $o$-diagonalizable.
A detailed study of $\fo$ is initiated in \cite{MShT:858},
where it is shown that consistently $\fo<\min\{\h,\s,\b\}$.

\begin{prob}[{Mildenberger, Shelah, et.\ al.\ \cite{MShT:858}}]
Is it consistent (relative to ZFC) that $\cov(\M)<\fo$?
\end{prob}

This problem, which originated from the topological
studies of the minimal tower problem, is of similar flavor: It is well-known that
if $\p=\aleph_1$, then $\t=\aleph_1$ too. We have a similar assertion for
$\cov(\M)$ and $\fo$: If $\cov(\M)=\aleph_1$, then $\fo=\aleph_1$, either.

\subsection{A table of open problems}
It is possible that the diagram in Figure \ref{tauSch} is incomplete:
There are many unsettled possible implications in it.
After \cite{tautau, ShTb768}, there remained $76$ (!) potential implications
which were not proved or ruled out.
The study made recently in \cite{MShT:858} and described in the previous section
ruled out $21$ of these implications, so that $55$ implications remain
unsettled.
The situation is summarized in Table \ref{imptab}, which updates the corresponding
table given in \cite{futurespm}.

Each entry $(i,j)$ ($i$th row, $j$th column) contains a symbol.
\checkmark means that property $(i)$ in Figure \ref{tauSch} implies
property $(j)$ in Figure \ref{tauSch}.
$\times$ means that property $(i)$ does not (provably) imply property $(j)$,
and \textbf{?} means that the corresponding implication is still unsettled.

\begin{table}[!ht]
\begin{changemargin}{-3cm}{-3cm}
\begin{center}
{\tiny
\begin{tabular}{|r||cccccccccccccccccccccc|}
\hline
   & \smb{0} & \smb{1} & \smb{2} & \smb{3} & \smb{4} & \smb{5} & \smb{6} & \smb{7} &
   \smb{8} & \smb{9} & \smb{10} & \smb{11} & \smb{12} & \smb{13} & \smb{14} & \smb{15} &
   \smb{16} & \smb{17} & \smb{18} & \smb{19} & \smb{20} & \smb{21}\cr
\hline\hline

\mb{ 0} &
\yup&\yup&\yup&\yup&\nop&\nop&\nop&\nop&\nop&\nop&\nop&
\nop&\yup&\yup&\nop&\mbq&\nop&\nop&\yup&\yup&\yup&\yup\\
\mb{ 1} &
\mbq&\yup&\yup&\yup&\nop&\nop&\nop&\nop&\nop&\nop&\nop&
\nop&\yup&\yup&\nop&\mbq&\nop&\nop&\mbq&\yup&\yup&\yup\\
\mb{ 2} &
\nop&\nop&\yup&\yup&\nop&\nop&\nop&\nop&\nop&\nop&\nop&
\nop&\nop&\yup&\nop&\mbq&\nop&\nop&\nop&\nop&\yup&\yup\\
\mb{ 3} &
\nop&\nop&\nop&\yup&\nop&\nop&\nop&\nop&\nop&\nop&\nop&
\nop&\nop&\nop&\nop&\nop&\nop&\nop&\nop&\nop&\nop&\yup\\
\mb{ 4} &
\yup&\yup&\yup&\yup&\yup&\yup&\yup&\yup&\nop&\nop&\mbq&
\mbq&\yup&\yup&\yup&\yup&\nop&\mbq&\yup&\yup&\yup&\yup\\
\mb{ 5} &
\mbq&\yup&\yup&\yup&\mbq&\yup&\yup&\yup&\nop&\nop&\mbq&
\mbq&\yup&\yup&\yup&\yup&\nop&\mbq&\mbq&\yup&\yup&\yup\\
\mb{ 6} &
\nop&\nop&\yup&\yup&\nop&\nop&\yup&\yup&\nop&\nop&\mbq&
\mbq&\nop&\yup&\nop&\yup&\nop&\mbq&\nop&\nop&\yup&\yup\\
\mb{ 7} &
\nop&\nop&\nop&\yup&\nop&\nop&\nop&\yup&\nop&\nop&\nop&
\mbq&\nop&\nop&\nop&\nop&\nop&\nop&\nop&\nop&\nop&\yup\\
\mb{ 8} &
\yup&\yup&\yup&\yup&\yup&\yup&\yup&\yup&\yup&\yup&\yup&
\yup&\yup&\yup&\yup&\yup&\yup&\yup&\yup&\yup&\yup&\yup\\
\mb{ 9} &
\mbq&\yup&\yup&\yup&\mbq&\yup&\yup&\yup&\mbq&\yup&\yup&
\yup&\yup&\yup&\yup&\yup&\yup&\yup&\mbq&\yup&\yup&\yup\\
\mb{10} &
\nop&\nop&\yup&\yup&\nop&\nop&\yup&\yup&\nop&\nop&\yup&
\yup&\nop&\yup&\nop&\yup&\nop&\yup&\nop&\nop&\yup&\yup\\
\mb{11} &
\nop&\nop&\nop&\yup&\nop&\nop&\nop&\yup&\nop&\nop&\nop&
\yup&\nop&\nop&\nop&\nop&\nop&\nop&\nop&\nop&\nop&\yup\\
\mb{12} &
\mbq&\mbq&\mbq&\mbq&\nop&\nop&\nop&\nop&\nop&\nop&\nop&
\nop&\yup&\yup&\nop&\mbq&\nop&\nop&\mbq&\yup&\yup&\yup\\
\mb{13} &
\nop&\nop&\nop&\nop&\nop&\nop&\nop&\nop&\nop&\nop&\nop&
\nop&\nop&\yup&\nop&\mbq&\nop&\nop&\nop&\nop&\yup&\yup\\
\mb{14} &
\mbq&\mbq&\mbq&\mbq&\nop&\nop&\nop&\nop&\nop&\nop&\nop&
\nop&\yup&\yup&\yup&\yup&\nop&\mbq&\mbq&\yup&\yup&\yup\\
\mb{15} &
\nop&\nop&\nop&\nop&\nop&\nop&\nop&\nop&\nop&\nop&\nop&
\nop&\nop&\yup&\nop&\yup&\nop&\mbq&\nop&\nop&\yup&\yup\\
\mb{16} &
\mbq&\mbq&\mbq&\mbq&\mbq&\mbq&\mbq&\mbq&\mbq&\mbq&\mbq&
\mbq&\yup&\yup&\yup&\yup&\yup&\yup&\mbq&\yup&\yup&\yup\\
\mb{17} &
\nop&\nop&\nop&\nop&\nop&\nop&\nop&\nop&\nop&\nop&\nop&
\nop&\nop&\yup&\nop&\yup&\nop&\yup&\nop&\nop&\yup&\yup\\
\mb{18} &
\nop&\nop&\nop&\nop&\nop&\nop&\nop&\nop&\nop&\nop&\nop&
\nop&\nop&\mbq&\nop&\mbq&\nop&\nop&\yup&\yup&\yup&\yup\\
\mb{19} &
\nop&\nop&\nop&\nop&\nop&\nop&\nop&\nop&\nop&\nop&\nop&
\nop&\nop&\mbq&\nop&\mbq&\nop&\nop&\nop&\yup&\yup&\yup\\
\mb{20} &
\nop&\nop&\nop&\nop&\nop&\nop&\nop&\nop&\nop&\nop&\nop&
\nop&\nop&\mbq&\nop&\mbq&\nop&\nop&\nop&\nop&\yup&\yup\\
\mb{21} &
\nop&\nop&\nop&\nop&\nop&\nop&\nop&\nop&\nop&\nop&\nop&
\nop&\nop&\nop&\nop&\nop&\nop&\nop&\nop&\nop&\nop&\yup\\

\hline
\end{tabular}
}
\end{center}
\end{changemargin}
\caption{Known implications and nonimplications}\label{imptab}
\end{table}

\begin{prob}[{Mildenberger, Shelah, et.\ al.\ \cite{tautau, MShT:858}}]
Settle any of the remaining $55$ implication in Table \ref{imptab}.
\end{prob}

\subsection{The Minimal Tower Problem revisited}
The study of $\tau$-covers was motivated by a combinatorial problem.
Interestingly, it rewarded us with a closely related \emph{purely combinatorial}
problem.
One of the difficulties with treating $\tau$-covers is the following:
Unlike the case of $\w$-covers or $\gamma$-covers, it could be that
$\cU$ is a $\tau$-cover of $X$ which refines another cover $\cV$ of $X$,
but $\cV$ is not a $\tau$-cover of $X$. This led to the introduction of
the following close relative \cite{tautau}.

A family $Y\sbst\roth$ is \Emph{linearly refinable}
if for each $y\in Y$ there exists an infinite subset
$\hat y\sbst y$ such that the family $\hat Y = \{\hat y : y\in Y\}$ is
linearly (quasi)ordered by $\as$.
A cover $\cU$ of $X$ is a \Emph{$\tau^*$-cover} of $X$ if it is large,
and $h[X]$ (where $h$ is the Marczewski characteristic function of $\cU$,
see Section \ref{dictionary}) is linearly refinable.
Let $\Tau^*$\index{$\Tau^*$} denote the collection of all \emph{countable} open
$\tau^*$-covers of $X$.

If we restrict attention to countable covers (this does not make a difference
for sets of reals), then $\Tau\sbst\Tau^*\sbst\Omega$.
Let $\p^*=\non(\binom{\Omega}{\Tau^*})$\index{$\p^*$}.
Then by the dictionary (Section \ref{dictionary}),
$$\p^* = \min\{|Y| : Y\sbst\roth\mbox{ is centered but not linearly refineable}\},$$
and we have the following interesting theorem (exercise):
\begin{thm}\label{pt}
$\p=\min\{\p^*,\t\}$.
\end{thm}
Thus $\p<\t$ implies the somewhat more pathological situation $\p^*<\t$.
If $\p^*$ is not provably equal to $\p$ then Theorem \ref{pt}
may be regarded as a partial (but easy) solution to the Minimal Tower Problem.
\begin{prob}[{Shelah, et.\ al.\ \cite{tautau, ShTb768}}]
Is $\p=\p^*$?
\end{prob}
As we have mentioned before, the closely related problem
whether $\p=\non(\binom{\Omega}{\Tau})$ has a positive answer
in \cite{ShTb768}.

\section{Some connections with other fields}

\subsection{Ramsey theory}\label{ramsey}

Recall that \Emph{Ramsey's Theorem}, often written as
$$\(\forall n,k\)\ \aleph_0\to\(\aleph_0\)^n_k,$$
asserts that for each $n$, $k$, and a countably infinite set $I$:
For each coloring $f:[I]^n\to\{1,\ldots,k\}$,
there exists a color $j$ and an infinite $J\sbst I$ such that $f\|[J]^n\equiv j$.

This motivates the following prototype for Ramsey theoretic hypotheses \cite{coc1}.
\bi
\item[$\scrA\to(\scrB)^n_k$\index{$\scrA\to(\scrB)^n_k$}:] For each $\cU\in\scrA$ and $f:[\cU]^n\to\{1,\ldots,k\}$,
there exists $j$ and $\cV\sbst\cU$ such that $\cV\in\scrB$ and $f\|[\cV]^n\equiv j$.
\ei
Using this notation, Ramsey's Theorem is $(\forall n,k)\ \roth\to(\roth)^n_k$.
This prototype can be applied to the case where $\scrA$ and $\scrB$ are collections
of rich covers, in accordance to the major theme of Ramsey theory that
often, when rich enough structures are split into finitely many pieces,
one of the pieces contains a rich substructure (Furstenberg).

The simplest case of the mentioned Ramseyan property is where $n=1$ (so that
we color the elements of the given member of $\scrA$ rather than finite subsets of it).
This case is well understood:
Fix a space $X$.
Since a finite partition of an infinite set must contain an infinite element,
and since every infinite subset of a $\gamma$-cover of $X$ is again a $\gamma$-cover of $X$,
we have that
$$(\forall k)\ \Gamma\to(\Gamma)^1_k.$$
It is less trivial but still not difficult to show that
$$(\forall k)\ \Omega\to(\Omega)^1_k$$
also holds, and since each $\tau$-cover is an $\omega$-cover (which in turn is a large cover),
and each large subcover of a $\tau$-cover is again a $\tau$-cover, we have that
$$(\forall k)\ \Tau\to(\Tau)^1_k.$$
The corresponding assertions for Borel or clopen (instead of open) covers also hold
for the same reasons, and it can be shown that nothing more can be said (except for
what trivially follows from the above assertions), concerning the partition property
$\scrA\to(\scrB)^1_k$ where $\scrA,\scrB\in\{\O,\Omega,\Tau,\Gamma\}$.
We now turn to the case that $n>1$.

Ramsey's Theorem is exactly the same as
$$\(\forall n,k\)\ \Gamma\to\(\Gamma\)^n_k.$$
However, if we substitute other classes of covers for $\scrA$ and $\scrB$
the property may not hold for all $X$.
An elegant example is the following.
\begin{thm}[{Just-Miller-Scheepers-Szeptycki \cite{coc1, coc2}}]\label{covRamsey}
The following properties are equivalent:
\be
\item $\sone(\Omega,\Omega)$,
\item $\Omega\to(\Omega)^2_2$; and
\item $(\forall n,k)\ \Omega\to(\Omega)^n_k$.
\ee
\end{thm}
In other words, $X$ satisfies $\sone(\Omega,\Omega)$ if, and only if,
whenever we color with $2$ colors the edges of a complete graph
whose vertices are elements of an $\omega$-cover of $X$,
we can find a complete monochromatic subgraph\index{complete monochromatic subgraph}
whose vertices form an $\omega$-cover of $X$.

\begin{cor}
~\be
\item $\sone(\Omega,\Gamma)$ is equivalent to $\Omega\to(\Gamma)^2_2$.
\item $\sone(\Omega,\Tau)$ implies $\Omega\to(\Tau)^2_2$.
\ee
And in both cases we can use any $n$ and $k$ as a superscript and a subscript.
\end{cor}
\begin{proof}
$(\Impl)$ for (1) and (2): Assume that $\Omega\spst\scrB$ and $X$ satisfies $\sone(\Omega,\scrB)$.
Then $X$ satisfies $\sone(\Omega,\Omega)$ as well as $\binom{\Omega}{\scrB}$.
By Theorem \ref{covRamsey}, $X$ satisfies $\Omega\to(\Omega)^2_2$,
but since $X$ satisfies $\binom{\Omega}{\scrB}$ and a (complete) subgraph of a
monochromatic graph is monochromatic, $X$ satisfies $\Omega\to(\scrB)^2_2$.

$(\Leftarrow)$ for (1): Clearly, $\Omega\to(\Gamma)^2_2$ implies $\binom{\Omega}{\Gamma}$,
which we know is the same as $\sone(\Omega,\Gamma)$.
\end{proof}

\begin{prob}
Is $\sone(\Omega,\Tau)$ equivalent to $\Omega\to(\Tau)^2_2$ ?
\end{prob}

Even for $n=k=2$, Ramsey's Theorem cannot be extended to $\aleph_1$.
To see this, consider a set of reals $X$ of size $\aleph_1$, and
a wellordering $\prec$ of $X$. Give a pair $\{x,y\}\sbst X$ the
color $1$ if $<$ and $\prec$ agree on $(x,y)$, and $0$ otherwise.
Then a monochromatic subgraph would give rise to an $<$-chain of reals
of size $\aleph_1$, which is impossible, since between each two
elements of the chain there is a rational number.
(This nice argument is due to Sierpi\'nski.)

There are some other ways to extend Ramsey's Theorem to larger cardinals,
sometimes taking a quick route to large cardinals, see \cite{HajLar}
for a survey. But Theorem \ref{covRamsey} and its relatives suggest
other generalizations. If we forget about the topology and consider
arbitrary (but countable) covers of the given space (that is, set),
then Theorem \ref{covRamsey} implies the following.
For an infinite cardinal $\kappa$, let $\Omega_\kappa$\index{$\Omega_\kappa$} denote
the collection of all \emph{countable} $\omega$-covers of $\kappa$.

\begin{thm}[{Scheepers \cite{SchBCT}}]
If $\kappa<\cov(\M)$, then $(\forall n,k)\ \Omega_\kappa\to(\Omega_\kappa)^n_k$,
where $\Omega_\kappa$ is the collection of countable $\w$-covers of $\kappa$.
(Moreover, for each $n$ and $k$ the bound $\cov(\M)$ is tight.)
\end{thm}
This is so because $\non(\sone(\Omega,\Omega))=\cov(\M)$.
This is a nice and typical case where the results in the topological studies,
which were motivated by results from pure infinite combinatorics, often project
to new results in pure infinite combinatorics.

\begin{rem}
For a much more comprehensive survey of this subject, see
\cite{KocRamsey}.
\end{rem}

\subsection{Countably distinct representatives and splittability}

Among the main tools for proving Ramsey theoretic results of the flavor of the
previous section one can find the properties of the following form:
\bi
\item[$\CDR(\scrA,\scrB)$\index{$\CDR(\scrA,\scrB)$}:] For each sequence $\seq{\cU_n}$ of elements of $\scrA$,
there exist countable, pairwise disjoint elements $\cV_n\sbst\cU_n$, $n\in\N$,
such that $\cV_n\in\scrB$ for all $n$.
\item[$\split(\scrA,\scrB)$\index{$\split(\scrA,\scrB)$}:] For each $\cU\in\scrA$, there exist disjoint $\cV_1,\cV_2\sbst\cU$
such that $\cV_1,\cV_2\in\scrB$.
\ei
Clearly, $\CDR(\scrA,\scrB)$ implies $\split(\scrA,\scrB)$,
of which an almost complete classification was carried in
\cite{splittability},
when $\scrA,\scrB\in\{\Omega,\Tau,\Gamma,\Lambda\}$.
As before, some of the properties are trivial, and several
equivalences are provable among the remaining properties.

The following dictionary is the key behind the classification,
where the negation of each property corresponds to some combinatorial
property of $h[X]$ where $h$ is the Marczewski characteristic function
(Section \ref{dictionary}) of $\cU$ for a cover $\cU$ witnessing the
failure of that property.

\begin{center}
\begin{tabular}{|c|c|}
\hline
Property & $h[X]$\\
\hline\hline
$\lnot\split(\Lambda,\Lambda)$ & reaping family\\
$\lnot\split(\Omega,\Lambda)$ & ultrafilter base\\
$\lnot\split(\Omega,\Omega)$ & ultrafilter subbase\\
$\lnot\split(\Tau,\Tau)$ & simple $P$-point base\\
\hline
\end{tabular}
\end{center}

The results of the classification are summarized in the
following theorem.

\begin{thm}[\cite{splittability}]\label{spl}
No additional implication (except perhaps the dotted ones) is provable.
$$\xymatrix{
\split(\Lambda, \Lambda) \ar[r] & \split(\Omega, \Lambda) \ar[r] & \split(\Tau, \Tau)\\
                           & \split(\Omega, \Omega)\ar[u]\\
& \split(\Omega, \Tau)\ar[u]\ar@{.>}[dr]^{(1)}\ar@/_/@{.>}[dl]_{(2)}\ar@{.>}[uul]_{(3)}\\
\split(\Omega, \Gamma) \ar[uuu]\ar[ur]\ar[rr]     & & \split(\Tau,\Gamma)\ar[uuu]\\
}$$
\end{thm}

\begin{prob}[\cite{splittability}]
Is the dotted implication (1) (and therefore (2) and (3)) in the diagram true?
If not, then is the dotted implication (3) true?
\end{prob}

\subsection{An additivity problem}
With regards to the additivity (preservation under taking finite unions) and
$\sigma$-additivity (countable unions), the following is known (\checkmark means that,
in the figure of Theorem \ref{spl},
the property in this position is $\sigma$-additive, and $\x$ means that it is not additive).
$$\xymatrix@R=10pt{
\txt{?} \ar[r] & \txt{\checkmark} \ar[r] & \txt{\checkmark}\\
                           & \x\ar[u]\\
& \x\ar[u]\\
\x \ar[uuu]\ar[ur]\ar[rr]     & & \txt{\checkmark}\ar[uuu]\\
}$$
Thus, the only unsettled problem is the following:
\begin{prob}[\cite{splittability}]\label{splitadd}
Is $\split(\Lambda,\Lambda)$ additive?
\end{prob}
In Proposition 1.1 of \cite{splittability} it is shown that
for a set of reals $X$ (in fact, for any hereditarily Lindel\"of space $X$),
each large open cover of $X$ contains a \emph{countable} large open cover of $X$.
Consequently, using standard arguments \cite{splittability}, the problem is closely
related to the following one (where $\roth$ is the space of all infinite
sets of natural numbers, with the topology inherited from $P(\N)$, the latter identified with $\Cantor$).
\begin{prob}
If $\mathsf{R}$ denotes the sets of reals $X$ such that
each continuous image of $X$ in $\roth$ is not reaping,
then is $\mathsf{R}$ additive?
\end{prob}

Zdomskyy has proved the following surprising result concerning the Hurewicz and Menger properties.

\begin{thm}[{Zdomskyy \cite{MenAdd}}]
Assume that $\u<\g$. Then $\ufin(\O,\Gamma)=\split(\Lambda,\Lambda)$, and $\ufin(\O,\Omega)=\ufin(\O,\O)$.
\end{thm}

Since $\ufin(\O,\Gamma)$ is easily seen to be countably additive, it follows that the answer
to Problem \ref{splitadd} is consistently positive.

\subsection{Topological games}\index{Topological games}
Historically, the bridge from general to infi\-nite-combinatorial topology
was through topological games related to covering properties,
introduced by Telg\'arsky in \cite{TelGames1, TelGames2},
and extensively studied by him and his colleagues (see Telg\'arsky's survey \cite{TelGamesSurv}).
This sort of game theory is still an important tool in proving Ramsey-theoretic results as those
in Section \ref{ramsey}.
The games appear under various guises, but we will focus on their form which is motivated by the
selection hypotheses.

$\gone(\scrA,\scrB)$\index{$\gone(\scrA,\scrB)$} is the game-theoretic version of $\sone(\scrA,\scrB)$.
In this game, ONE chooses in the $n$th inning an element
$\cU_n\in\scrA$ and then TWO responds by choosing $U_n\in \cU_n$.
They play an inning per natural number.
This is illustrated in the following figure.

\medskip

\begin{center}
\begin{tabular}{cccccccccc}
ONE: & $\cU_1\in\scrA$ &    &      &      & $\cU_2\in\scrA$ &     & & \dots\\
     &     & $\searrow$    &     & $\nearrow$    &      & $\searrow$\\
TWO: &   &    & $U_1\in\cU_1$ &   &     &    & $U_2\in\cU_2$ & \dots\\
\end{tabular}
\end{center}

\medskip

\noindent TWO wins if $\{U_1,U_2,\dots\}\in\scrB$, otherwise ONE wins.

The game $\gfin(\scrA,\scrB)$\index{$\gfin(\scrA,\scrB)$} is played similarly, where TWO responds with
finite subsets $\cF_n\sbst\cU_n$ and wins if $\Union_n\cF_n\in\scrB$.

Observe that if ONE does not have a winning strategy in
$\gone(\scrA,\scrB)$ (respectively, $\gfin(\scrA,\scrB)$), then $\sone(\scrA,\scrB)$
(respectively, $\sfin(\scrA,\scrB)$) holds.
The converse is not always true; when it is true,
the game is a powerful tool for studying the combinatorial
properties of $\scrA$ and $\scrB$.
Fortunately, this is often the case.
In fact, this is always the case for the properties in the Scheepers Diagram
(see, e.g., \cite{coc7, coc8, strongdiags}, and references therein).

There is, though, a well known property of similar flavor for which
the question of equivalence is still open.
Consider the following generalized selection hypothesis.
\begin{itemize}
\item[$\sone(\seq{\scrA_n},\scrB)$:]\index{$\sone(\seq{\scrA_n},\scrB)$}
For each sequence of elements $\cU_n\in\scrA_n$, $n\in\N$,
there are elements $U_n\in\cU_n$, $n\in\N$, such that $\{U_n : n\in\N\}\in\scrB$.
\end{itemize}
Define its game theoretic version $\gone(\seq{\scrA_n},\scrB)$\index{$\gone(\seq{\scrA_n},\scrB)$}
in the natural way.
For most instances of $\seq{\scrA_n}$ and $\scrB$, we do not get anything new by considering
these properties \cite{strongdiags}. However, this is not always the case.

A cover $\cU$ of $X$ is an \emph{$n$-cover} of $X$ if each subset $F$ of $X$
of cardinality at most $n$ is contained in some member of $\cU$.
For each $n$ denote by $\O_n$ the collection of all
open $n$-covers of $X$.
Then the property $\sone(\seq{\O_n},\Gamma)$\index{strong $\gamma$-property}, introduced
by Galvin and Miller in \cite{GM} (where it was called \emph{strong $\gamma$-property}),\footnote{Actually,
Galvin and Miller defined the property as $(\exists k_n\nearrow\infty)\ X\in\sone(\seq{\O_{k_n}},\Gamma)$,
but it was observed in \cite{strongdiags} that the quantifier can be eliminated from their
definition.
}
is strictly stronger than
$\sone(\Omega,\Gamma)$
(which in turn is the strongest property in the Scheepers Diagram) \cite{BR}.
The strong $\gamma$-property $\sone(\seq{\O_n},\Gamma)$ never had a game theoretic
characterization. A natural candidate is the following.
Let $\gone(\seq{\scrA_n},\scrB)$\index{$\gone(\seq{\scrA_n},\scrB)$}
be the game theoretic version of $\sone(\seq{\scrA_n},\scrB)$
(so it is like $\gone(\scrA,\scrB)$, but in the $n$th inning ONE chooses $\cU_n\in\scrA_n$
instead of $\cU_n\in\scrA$).

\begin{prob}[\cite{strongdiags}]
Is the strong $\gamma$-property $\sone(\seq{\O_n},\Gamma)$  equivalent to ONE not having
a winning strategy in $\gone(\seq{\O_n},\Gamma)$ ?
\end{prob}

What about the case that TWO has a winning strategy in the game
$\gone(\scrA,\scrB)$ or $\gfin(\scrA,\scrB)$?
It turns out that in the cases corresponding to the properties in the Scheepers Diagram,
these questions have interesting and elegant solutions.
In this interpretation, the conjectures made by Menger, Hurewicz, and
Borel are all correct!\index{Borel Conjecture}\index{Menger Conjecture}\index{Hurewicz Conjecture}

\begin{thm}[Telg\'arsky]
For a metrizable space $X$:
TWO has a winning strategy in the game $\gfin(\O,\O)$ if, and only if,
the space $X$ is $\sigma$-compact.
\end{thm}
Note that $\sfin(\O,\O)$ is the same as the Menger property $\ufin(\O,\O)$,
so $\gfin(\O,\O)$ is the game theoretic version of the Menger property.

\begin{thm}[Galvin; Telg\'arsky]
For a metrizable space $X$:
TWO has a winning strategy in the game $\gone(\O,\O)$ if, and only if,
the space $X$ is countable.
\end{thm}
$\gone(\O,\O)$ is the game theoretic version of Rothberger' property $\sone(\O,\allowbreak\O)$
\index{Rothberger property}. Recall from Section \ref{BCrevisited},
that the Borel Conjecture is equivalent to
the assertion that all spaces satisfying $\sone(\O,\O)$ are countable.

Another way to interpret these results is as follows:
If we assume that $\gone(\O,\O)$ is determined,
then Borel's Conjecture is true.

\begin{cor}
$\gone(\O,\O)$ is determined if, and only if, Borel's Conjecture is true.
\end{cor}

Note that, by Section \ref{ZFCexamples}, we have to give up the Axiom of Choice
in order to make the corresponding assertion for $\gfin(\O,\O)$ and the Menger Conjecture
meaningful.
Scheepers told us that the Axiom of Determinacy (which rules out the Axiom of Choice)
implies that the above-mentioned games are determined.

\subsection{Arkhangel'ski\v{i} duality theory}\label{fspaces}

Consider $C(X)$\index{$C(X)$}, the space of continuous real-valued functions
$f:X\to\R$, as a subspace of $\R^X$\index{$\R^X$} (the Tychonoff product of $X$
many copies of $\R$). This is often called the
\Emph{topology of pointwise convergence} and sometimes denoted $C_p(X)$\index{$C_p(X)$}
rather than $C(X)$.

The object $C(X)$ is very complicated from a topological point of view,
in fact, even the behavior of the closure operator in this space is
complicated. To this end, a comprehensive duality theory was developed
by Arkhangel'ski\v{i} and his followers, which translates properties
of $C(X)$ into properties of $X$, which are easier to work with.

For example, recall that a space $Y$ is \Emph{Fr\'echet-Urysohn} if
for each $A\sbst Y$ and $x\in\cl{A}$, there is a sequence $\seq{a_n}$ in $A$
such that $\lim_{n\to\infty} a_n=x$.
In their celebrated 1982 paper \cite{GN},
Gerlits and Nagy proved that $C(X)$ is Fr\'echet-Urysohn if, and only if,
$X$ satisfies $\binom{\Omega}{\Gamma}$ (or, equivalently, $\sone(\Omega,\Gamma)$).

An interesting example of the applicability of infinite-combinatorial methods for
these questions is the following.
In 1992, at a seminar in Moscow, Reznicenko introduced the following property:
A space $Y$ is \Emph{weakly Fr\'echet-Urysohn} if,
for each $A\sbst Y$ and $x\in\cl{A}\sm A$, there exist
finite disjoint sets $F_n\sbst A$, $n\in\N$,
such that for each neighborhood $U$ of $x$,
$F_n\cap U\neq\emptyset$ for all but finitely many $n$.

In \cite{FunRez}, Ko\v{c}inac and Scheepers made the following conjecture,
which is now a theorem.
\begin{thm}[\cite{reznicb}]
The minimal cardinality of a set of reals $X$ such that $C(X)$ does not have
the weak Fr\'echet-Urysohn property is $\b$.
\end{thm}
The surprising thing is that its proof only required a translation into the language
of combinatorics and an application of an existing result:
It was known that this minimal cardinality is at least $\b$.
A result of Sakai \cite{SakaiRez} can be used to prove that
if $C(X)$ is weakly Fr\'echet-Urysohn, then a continuous image of
$X$ cannot be a subbase for a non-feeble filter on $\N$
(see \cite{BlassHBK} for the definition of non-feeble filter).
Now it remained to apply a result of Petr Simon which tells that
there exists a non-feeble filter with base of size $\b$.

\section{Conclusions}
The theory emerging from the systematic study of diagonalizations
of covers in a unified framework is not only aesthetically pleasing,
but is also useful in turning otherwise ingenious ad hoc arguments
into natural explanations.
For this a good terminology is required, and the major part of this
was already suggested by Scheepers' selection prototypes.

This study has connections and applications in several related
fields, like Ramsey theory, function spaces, and topological groups.
The usage of infinite-combinatorial methods in this theory has proved
successful, and is often the ``correct'' tool to investigate these problems.
This approach sometimes rewords
by implying interesting results in the field of pure infinite-combinatorics.

While with regards to some of the investigations concerning the
classical types of covers the picture is rather complete now,
there remains much to be explored with regards to the new types of
covers and their connection to the related fields.

We have only given a tiny sample of each theme.
The reader is referred to \cite{LecceSurvey, KocSurv, futurespm} to have
a more complete picture of the framework and the open problems
it poses.

\subsection{Acknowledgments}
We thank Tomasz Weiss for the proof of Proposition \ref{weiss}.
We also thank Rastislav Telg\'arsky, Lyubomyr Zdomskyy, and the referee,
for making several interesting comments.

\ed